\documentclass{amsart}
\usepackage{array}
\usepackage{multirow}
\usepackage{graphicx}
\usepackage{epsfig}
\usepackage{amsfonts,amsmath}
\usepackage{amssymb}
\usepackage{latexsym}
\usepackage{epsf}
\usepackage{graphicx,color,graphics}
\usepackage{amsmath,amssymb}
\usepackage{dsfont}
\usepackage[latin1]{inputenc}  
\newtheorem{theorem}{Theorem}[section]

\newtheorem{lemma}[theorem]{Lemma}

\newtheorem{remark}{Remark}

\numberwithin{equation}{section}

\newcommand{\R}{\mathbb{R}}

\begin{document}

\title{Laminated Timoshenko beams with interfacial slip and infinite memories}

\pagenumbering{arabic}

\maketitle

\begin{center}
\vspace{0.3cm}
A. Guesmia\\
Institut Elie Cartan de Lorraine, UMR 7502, Universit\'e de Lorraine, 3 Rue Augustin Fresnel\\ 
BP 45112, 57073 Metz Cedex 03, France\\
and\\
Department of Mathematics and Statistics, King Fahd University of Petroleum and Minerals\\
Dhahran 31261, Saudi Arabia\\ 
e-mail: aissa.guesmia@univ-lorraine.fr\\
\vspace{0.3cm}
J. E. Mu\~noz Rivera\\
Department of Mathematics, Federal University of Rio de Janeiro and National Laboratory\\ 
for Scientific Computation, Petr\'opolis, RJ, Brasil\\
e-mail: rivera@lncc.br\\
\vspace{0.3cm}
M. A. Sep\'ulveda Cort\'es\\
CI$^2$MA and DIM, Universidad del Concepci\'{o}n, Concepci\'{o}n, Chile\\
e-mail: mauricio@ing-mat.udec.cl \\
\vspace{0.3cm}
O. Vera Villagr\a'{a}n\\
Departamento de Matem\a'{a}tica, Universidad del B\'{\i}o-B\'{\i}o, Concepci\'{o}n, Chile\\
e-mail: overa@ubiobio.cl
\end{center}
\begin{abstract}
We study in this paper the well-posedness and stability of three structures with interfacial slip and two infinite memories effective on the transverse displacement and the rotation angle. We consider a large class of kernels and prove that the system has a unique solution satisfying some regularity properties. Moreover, without restrictions on the values of the parameters, we show that the solution goes to zero at infinity and give an information on its speed of convergence in terms of the growth of kernels at infinity. 
A numerical analysis of the obtained theoretical results will be also given.
\end{abstract}

{\it keywords.} Viscoelastic structure, interfacial slip, semigroups, stability, lyapunov functional, numerical analysis.


\section{Introduction}

We consider the well-posedness and stability of three structures with interfacial slip and two infinite memories effective on the transverse displacement and the rotation angle 
\begin{align}\label{109}
\begin{cases}
\rho_{1}\,\varphi_{tt} + k\,(u - \varphi_{x})_{x} + \displaystyle\int_{0}^{\infty}\,g_{1}(s)\,\varphi_{xx} (t - s)\ ds = 0,  \\
\rho_{2}\,(v - u)_{tt} - b\,(v - u)_{xx} - k\,(u - \varphi_{x})  + \displaystyle\int_0^{\infty}\,g_{2} (s)\,(v(t - s) - u(t - s))_{xx}\ ds =0, \\
\rho_{2}\,v_{tt} - b\,v_{xx} + 3\,k\,(u - \varphi_{x}) + 4\,\delta\,v + 4\,\gamma\,v_{t}= 0
\end{cases}
\end{align}
with boundary conditions
\begin{equation}
\label{114}\varphi(0,\,t) = \varphi (1,\,t)= u (0,\,t)=u_{x} (1,\,t)=v(0,\,t)= v_{x} (1,\,t)=0
\end{equation}
and initial data
\begin{align}\label{115}
\begin{cases}
(\varphi(x,\,-t),\,u(x,\,-t),\,v(x,\,-t))=(\varphi_{0} (x,t),\,u_{0}(x,t),\,v_{0} (x,t)),\\
(\varphi_{t}(x,\,0),\,u_{t}(x,\,0),\,v_{t}(x,\,0))= (\varphi_{1}(x),\,u_{1} (x),\,v_{1} (x)),
\end{cases}
\end{align}
where $(x,\,t) \in (0,\,1) \times \mathbb{R}_{+}$, $\varphi =\varphi (x,\,t)$ is the transversal displacement, $u = u(x, \,t)$ represents the rotation angle, $v = v(x,\,t)$ is proportional to the amount of slip along the interface, $\rho_{1},$ $\rho_{2},$ $k,$ $b,$ $\delta$ and $\gamma$ are positive constants and $g_{i} :\mathbb{R}_{+}\to \mathbb{R}_{+}$ is a given function, $i=1,\,2$. 
\vskip0,1truecm
The structures known under the name laminated Timoshenko beams are composed of two layered identical beams of uniform thikness and attached together on top of each other subject to transversal and rotational vibrations, and taking account the longitudinal dislacement. An adhesive layer of small thikness is bonding the two adjoining surfaces and creating a restoring force being proportional to the amount of slip and producing a damping. These structures are used in many practical fields; see, for example, \cite{lota1} for more details. When the longitudinal dislacement (slip) is ignored, the laminated Timoshenko beams are reduced to the well known Timoshenko beams \cite{timo}. During the last few years, these structures were the subject of several studies in the literature recovering well-posedness and stability by adding some kinds of (boundary or internal) controls.
\vskip0,1truecm
When $v$ is not taken in consideration: $v=0$, \eqref{109} is reduced to Timoshenko beams \cite{timo} and its stability question was widely treated in a huge number of works; see, for example, \cite{cdflr} and \cite{gms} and the refereces therein. 
\vskip0,1truecm
The authors of \cite{wxy} proved the exponential stability through mixed homogeneous Dirichlet-Neumann boundary conditions and two boundary controls at $x=1$ provided that the speeds of wave propagations of the first two equations are different; that is
\begin{equation}
{\sqrt{\frac{k}{\rho_1}}}\ne {\sqrt{\frac{b}{\rho_2}}}. \label{wavespeeds}
\end{equation}
It was also proved in \cite{wxy} that the frictional damping $4\gamma v_t$ is strong enough to stabilize asymptotically the structure but it is not able to stabilize the structure exponentially. The same exponential stability result of \cite{wxy} was proved in \cite{clx} for the same model but through two boundary controls at $x=0$ and $x=1$. The expeonential stability result of \cite{wxy} was improved in \cite{tata} by assuming some weaker conditions on the parameters. The authors of \cite{lota3} proved that the exponential stability holds if the boundary controls are replaced by a frictional damping acting on the first equation. The author of \cite{rapo} proved that the exponential stability holds without any restriction on the parameters if the first two equations are also damped via frictional dampings. Recently, the authors of \cite{rvma} proved that, without any restriction on the parameters, the polynomial stability holds under additional three dynamic boundary conditions. For the stability of laminated beams with Cattaneo's or Fourier's type heat conduction, we refer the readers to \cite{agg} and \cite{lizh}.
\vskip0,1truecm
The stability in case of viscoelastic dampings represented by finite memory terms in the form of a convolution on $[0,t]$ (see, for example, \cite{beim}, \cite{lota1} and \cite{lota3}) was treated in \cite{lota1}, \cite{lota2} and \cite{lota3}. Namely, under some restrictions on the parameters and with one or two kernels converging exponentially to zero at infinity, the exponential stability was proved in \cite{lota1}, \cite{lota2} and \cite{lota3}.  
\vskip0,1truecm
For the stability of Bresse systems \cite{bres} with infinite memories, we refer the readers to \cite{gues2}, \cite{gues4}, \cite{gues6} and \cite{gues3}, and the references therein. 
\vskip0,1truecm
In \cite{gues8}, the first author of the present paper proved, under some restrictions on the parameters and with $\gamma=0$, some exponential and polynomial stability results using only one infinite memory with a kernel converging exponentially to zero at infinity.
\vskip0,1truecm
From the cited results above, we see that the exponential and/or polynomial stability has been proved under some restrictions on the parameters and with kernels converging exponentially to zero at infinity. The first main objective of this paper is proving that two infinite memories guarantee the stability of the system without any restrictions on the parameters. Moreover, we show that the class of admissible kernels is much more larger than the one containing kernels converging exponentially to zero at infinity. Our second main objective is presenting a numerical analysis to illustrate our theoretical stability results. 
\vskip0,1truecm
The proof of the well-posedness and stability results are based on the semigroup theory and the energy method, respectively. However, the numerical results are proved using the finite difference approach (of second order in space and time).
\vskip0,1truecm
The paper is organized as follows. In Section 2, we consider some assumptions on the relaxation functions and prove the well-posedness. In Sections 3, we prove our stability results. In Section 4, we present our numerical analysis. We end our paper by giving some general comments and issues in section 5.

\section{Setting of the Semigroup}

We introduce the variable $\psi =v-u$, and as in \cite{dafe}, we consider the variables $\eta$ and $z$ and theirs initial data given by
\begin{align*}
\begin{cases}
\eta (x,\,t,\,s)= \varphi (x,\,t) - \varphi (x,\,t - s),\quad & x\in (0,\,1),\ s,\,t\in \mathbb{R}_{+},    \\
\eta_{0}(x,\,s) =\varphi_{0}(x,\,0) - \varphi_{0} (x,\,s),\quad & x\in (0,\,1),\ s,\,t\in \mathbb{R}_{+},  \\
z (x,\,t,\,s) = \psi(x,\,t) - \psi (x,\,t - s), \quad & x\in (0,\,1),\ s,\,t\in \mathbb{R}_{+} \\ 
z_{0}(x,\,s) = v_{0}(x,\,0) - u_{0}(x,\,0) - (v_{0}(x,\,s) - u_{0}(x,\,s)),\quad & x\in (0,\,1),\ s,\,t\in \mathbb{R}_{+}.
\end{cases}
\end{align*}
So the system \eqref{109} becomes 
\begin{align}\label{109+}
\begin{cases}
\rho_{1}\,\varphi_{tt} - k\,(\varphi_{x} +\psi - v)_{x} + g_{1}^{0}\,\varphi_{xx} - \displaystyle\int_{0}^{\infty}\,g_{1}(s)\,\eta_{xx}(t - s)\ ds = 0, \\
\rho_{2}\,\psi_{tt} - (b - g_{2}^{0} )\,\psi_{xx} + k\,(\varphi_{x} + \psi - v)  - \displaystyle\int_{0}^{\infty}\,g_{2}(s)\,z_{xx}\ ds =0, \\
\rho_{2}\,v_{tt} - b\,v_{xx} - 3\,k\,(\varphi_{x} + \psi - v) + 4\,\delta\,v + 4\,\gamma\,v_{t} = 0,
\end{cases}
\end{align}
where 
\begin{equation*}
g_{i}^{0} =\int_{0}^{\infty}\,g_{i}(s)\ ds,\quad i=1,\,2,
\end{equation*}
with boundary conditions
\begin{equation}
\label{114+}\varphi (0,\,t) = \varphi (1,\,t) = \psi (0,\,t)=\psi_{x} (1,\,t) = v(0,\,t) =  v_{x} (1,\,t)=0,\quad t\in \mathbb{R}_{+}.
\end{equation}
The functionals $\eta$ and $z$ satisfy
\begin{align}\label{eta}
\begin{cases}
\eta_{t}(x,\,t,\,s) + \eta_{s}(x,\,t,\,s) - \varphi_{t}(x,\,t)=0,\quad & x\in (0,\,1),\  s,\,t>0,\\
z_{t}(x,\,t,\,s) + z_{s}(x,\,t,\,s) - \psi_{t}(x,\,t) = 0,\quad & x\in (0,\,1),\ s,\,t>0, \\
\eta(x,\,0,\,s) = \eta_{0}(x,\,s),\quad z (x,\,0,\,s) = z_{0}(x,\,s),\quad & x\in (0,\,1),\ s\in \mathbb{R}_{+}, \\
\eta (x,\,t,\,0) = z(x,\,t,\,0) = z_{x}(1,\,t,\,s) = 0,\quad & x\in (0,\,1),\, t,\,s\in \mathbb{R}_{+} .
\end{cases}
\end{align}
Let 
\begin{align}\label{U}
\begin{cases}
\varphi_{t} = \tilde{\varphi},\quad \psi_{t} =\tilde{\psi},\quad v_{t} = \tilde{v},\\
U = (\varphi,\,\tilde{\varphi},\,\psi,\,\tilde{\psi},\,v,\,\tilde{v},\,\eta,\,z),\\
U_0 =(\varphi_{0},\,\varphi_{1},\,v_{0} - u_{0},\,v_{1} - u_{1},\,v_{0},\,v_{1},\,\eta_{0},\,z_{0}).
\end{cases}
\end{align}
Now, we can rewrite the system \eqref{115}, \eqref{109+} and \eqref{114+} in the following initial value problem:
\begin{align}\label{222}
\begin{cases}
U_{t}(t) = \mathcal{A}U(t),\quad t>0, \\
U(0) = U_{0},
\end{cases}
\end{align}
where the operator $\mathcal{A}$ is defined by
\begin{equation}
\label{202}\mathcal{A}U =\left(
\begin{array}{c}
\tilde{\varphi} \\
\\
\frac{k}{\rho_{1}}\,(\varphi_{x} + \psi - v)_{x} - \frac{g_{1}^{0}}{\rho_{1}}\,\varphi_{xx} + \frac{1}{\rho_{1}}\int_{0}^{\infty}\,g_1 (s)\,\eta_{xx}\ ds \\
\\
\tilde{\psi} \\
\\
\frac{1}{\rho_{2}}\left[(b - g_{2}^{0} )\,\psi_{xx} - k\,(\varphi_{x} + \psi - v) \right] +\frac{1}{\rho_{2}}\int_{0}^{\infty}\,g_{2}(s)\,z_{xx}\ ds \\
\\
\tilde{v} \\
\\
\frac{1}{\rho_{2}}\left[b\,v_{xx} + 3\,k\,(\varphi_{x} + \psi -v) - 4\,\delta\,v - 4\,\gamma\,\tilde{v} \right]  \\
\\
-\,\eta_{s} + \tilde{\varphi} \\
\\
-\,z_{s} + \tilde{\psi} 
\end{array}
\right).
\end{equation}
\\
\\
Let us consider the standard $L^{2}(0,\,1)$ space with its classical scalar product $\langle\cdot,\,\cdot\rangle$ and generated norm $\Vert\cdot\Vert$. We consider also the phase Hilbert spaces
\begin{equation*}
L_{g_{i}} =\left\{v:\mathbb{R}_{+}\to {\tilde H}_{i},\ \int_{0}^{\infty}\,g_{i}(s)\,\Vert v_{x} (s)\Vert^{2}\ ds <\infty\right\},
\end{equation*}
equipped with the inner product
\begin{equation*}
\langle w,\,\tilde{w}\rangle_{L_{g_{i}}} = \int_{0}^{\infty}g_{i}(s)\,\langle {w}_{x} (s),\,\tilde{w}_{x} (s)\rangle\ ds, 
\end{equation*}
$\tilde{H}_{1} = H_{0}^{1}(0,\,1)\ $, $\tilde{H}_{2} = H_{*}^{1}(0,\,1)$ and
\begin{equation*}
H_{*}^{1}(0,\,1) = \{h\in H^{1}(0,\,1):\ h(0) = 0\}.
\end{equation*}
The energy space is given by 
\begin{equation*}
\mathcal{H}=H_{0}^{1}(0,\,1)\times L^{2}(0,\,1)\times[H_{*}^{1}(0,\,1) \times L^{2}(0,\,1)]^{2}\times L_{g_{1}} \times L_{g_{2}}
\end{equation*}
equipped with the inner product, for any 
\begin{equation*}
U_{1} =(\varphi_{1},\,\tilde{\varphi}_{1},\,\psi_{1},\,\tilde{\psi}_{1},\,v_{1},\,\tilde {v}_{1},\,\eta_{1},\,z_{1}),\quad U_{2} = (\varphi_{2},\,\tilde{\varphi}_{2},\,\psi_{2},\,\tilde{\psi}_{2},\,v_{2},\,\tilde{v}_{2},\,\eta_{2},\,z_{2})\in \mathcal{H}, 
\end{equation*}
\begin{eqnarray}
\label{303}
\langle U_{1} ,\,U_{2}\rangle_{\mathcal{H}} =  3\,k\,\langle (\varphi_{1x} +\psi_{1} - v_{1}),\,(\varphi_{2x} + \psi_{2} - v_{2})\rangle + 3\,(b - g_{2}^{0})\langle \psi_{1x},\,\psi_{1x}\rangle - 3\,g_{1}^{0}\langle \varphi_{1x},\,\varphi_{2x}\rangle 
\end{eqnarray}
\begin{equation*}
+\,b\,\langle v_{1x},\,v_{2x}\rangle + 4\,\delta\,\langle v_{1},\,v_{2}\rangle + 3\,\rho_{1}\langle \tilde{\varphi}_{1},\,\tilde{\varphi}_{2}\rangle + 3\,\rho_{2}\,\langle \tilde{\psi}_{1},\,\tilde{\psi}_{2}\rangle  + \rho_{2}\,\langle \tilde{v}_{1},\,\tilde{v}_{2} \rangle + 3\,\langle\eta_{1},\,\eta_{2}\rangle_{L_{g_{1}}} + 3\,\langle z_{1},\,z_{2}\rangle_{L_{g_{2}}}.
\end{equation*}
The domain of $\mathcal{D}(\mathcal{A})$ is given by
\begin{equation}
\label{3030}\mathcal{D}(\mathcal{A}) = \{U\in \mathcal{H},\, \mathcal{A}U\in \mathcal{H},\ \psi_{x}(1) = v_{x}(1)=\eta(x,\,0) = z(x,\,0) = z_{x}(1,\,s)=0\}.
\end{equation}
Now, to get our well-posedness results, we consider the following hypothesis: 
\vskip0,1truecm
${\bf (H1)}$ Assume that the function $g_{i}:\mathbb{R}_{+}\to \mathbb{R}_{+}$, $i=1,\,2$, is differentiable, nonincreasing and integrable on $\mathbb{R}_{+}$ such that there exists a positive constant $k_{0}$ satisfying, for any 
\begin{equation*}
(\varphi,\,\psi,\,w)\in H_{0}^{1}(0,\,1)\times H_{*}^{1}(0,\,1)\times H_{*}^{1} (0,\,1),
\end{equation*} 
we have
\begin{eqnarray}
\label{2.16}
k_{0}\left(\Vert\varphi_{x}\Vert^{2} + \Vert\psi_{x}\Vert^{2} + \Vert w_{x}\Vert^{2} \right) & \leq & 3\,k\,\Vert\varphi_{x} + \psi - v\Vert^{2} + 3\,(b - g_{2}^{0})\,\Vert\psi_{x}\Vert^{2} \nonumber \\
& & +\ b\,\Vert v_{x}\Vert^{2} + 4\,\delta\,\Vert v\Vert^{2} - 3\,g_{1}^{0}\,\Vert \varphi_{x}\Vert^{2} .
\end{eqnarray}
Moreover, assume that there exists positive constants $\beta_{1}$ and $\beta_{2}$ such that
\begin{equation}\label{2.230}
-\,\beta_{i}\,g_{i}(s) \leq g'_{i}(s),\quad s\in \mathbb{R}_{+} , \quad i=1,\,2.
\end{equation} 
\vskip0,1truecm
\begin{theorem}\label{Theorem 1.1}
Assume that ${\bf (H1)}$ holds. Then, for any $U_{0} \in \mathcal{D}(\mathcal{A})$, system \eqref{222} admits a unique solution $U$ satisfying
\begin{equation}
U\in C \left(\mathbb{R}_{+} ;\mathcal{H}\right). \label{exist1}
\end{equation}
If $U_{0} \in \mathcal{D}(\mathcal{A})$, then $U$ satisfies
\begin{equation}
U\in C^{1} \left(\mathbb{R}_{+};\,\mathcal{H}\right)\cap C \left(\mathbb{R}_{+};\,\mathcal{D}\left(\mathcal{A}\right)\right). \label{exist2}
\end{equation}
\end{theorem}
\vskip0,1truecm
\textbf{Proof}. First, under condition \eqref{2.16} and according to our homogeneous Dirichlet boundary conditions and Poincare\'e's inequality, the inner product of $\mathcal{H}$ generates a norm
\begin{equation}\label{303}
\Vert U\Vert_{\mathcal{H}}^{2} =  \langle U,\,U\rangle_{\mathcal{H}}
\end{equation} 
equivalent to the one of $[H^1 (0,\,1)\times L^{2} (0,\,1)]^{3}\times L_{g_{1}}\times L_{g_{2}}$, and so $\mathcal{H}$ is a Hilbert space.
The proof of Theorem \ref{Theorem 1.1} relies then on the Lumer-Philips theorem by proving that the operator $\mathcal{A}$ is dissipative and $I-\mathcal{A}$ is surjective ($I$ denotes the identity operator); that is $-\mathcal{A}$ is maximal monotone. So $\mathcal{A}$ is then the infinitesimal generator of a C$_{0}$-semigroup of contraction on $\mathcal{H}$ and its domain $\mathcal{D}(\mathcal{A})$ is dense in $\mathcal{H}$. The conclusion then follows immediately (see \cite{pazy}).
\vskip0,1truecm
Second, direct calculations give
\begin{equation}
\left\langle \mathcal{A}U,\,U \right\rangle_{\mathcal{H}} = -\,4\,\gamma\,\Vert v_{t}\Vert^{2} + \frac{1}{2}\int_{0}^{\infty}g'_{1}(s)\left\Vert \eta_{x}\right\Vert^{2} ds + \frac{1}{2}\int_{0}^{\infty}g'_{2} (s)\left\Vert z_{x}\right\Vert^{2}\ ds \leq 0. \label{dissp}
\end{equation}
Hence, $\mathcal{A}$ is a dissipative operator, since $g_{1}$ and $g_{2}$ are nonincreasing and \eqref{2.230} guarantees the boundedness of the integrals in \eqref{dissp}.
\vskip0,1truecm
Third, we prove that $I - \mathcal{A}$ is surjective. Let $F=(f_{1},\,\ldots,\,f_{8} )^{T} \in \mathcal{H}$. We prove that there exists 
$U \in \mathcal{D}\left( \mathcal{A}\right)$ satisfying
\begin{equation}
U - \mathcal{A}U = F.  \label{ZF}
\end{equation}
First, the first, third and fifth equations in \eqref{ZF} are equivalent to
\begin{equation}
\tilde{\varphi} = \varphi - f_{1} ,\quad \tilde{\psi} =\psi - f_{3} \quad\hbox{and}\quad \tilde{v} =v - f_{5}. \label{z1f1}
\end{equation}
Second, from \eqref{z1f1}, we see that the last two equations in \eqref{ZF} are reduced to 
\begin{equation}
\eta_{s} + \eta = \varphi + f_{7} - f_{1} \quad\hbox{and}\quad z_{s} + z = \psi + f_{8} - f_{3}. \label{z7*}
\end{equation}
Integrating with respect to $s$ and noting that $\eta$ and $z$ should satisfy $\eta (0) =z(0) = 0$, we get
\begin{equation}
\eta (s) = (1 - e^{-\,s})(\varphi -f_{1}) + \int_{0}^{s} e^{\tau - s}\,f_{7} (\tau)\ d\tau\,\,\hbox{and}\,\,
z(s) = (1 - e^{-\,s})(\psi -f_{3}) + \int_{0}^{s} e^{\tau - s}\,f_{8} (\tau)\ d\tau. \label{z7}
\end{equation}
Third, using \eqref{z1f1} and \eqref{z7}, we find that the second, fourth and sixth equations in $\eqref{ZF}$
are reduced to
\begin{equation}
\left\{
\begin{array}{ll}
\label{z5f5}
\rho_{1}\,\varphi - k\,(\varphi_{x} + \psi - v)_{x} + \tilde{g}_{1}\,\varphi_{xx} = \rho_{1} (f_{1} + f_{2}) +
\displaystyle\int_{0}^{\infty}g_{1}(s)\left((1 - e^{-\,s})f_{1} + \displaystyle\int_{0}^{s}e^{\tau - s}\,f_{7}(\tau)\,d\tau\right)_{xx}\ ds, \vspace{0.2cm}  \\
\rho_{2}\,\psi - (b - \tilde{g}_{2} )\psi_{xx} + k\,(\varphi_{x} + \psi - v) = \rho_{2}\, (f_{3} + f_{4}) +
\displaystyle\int_{0}^{\infty}g_{2}(s)\left((1 - e^{-\,s})\,f_{3} + \displaystyle\int_{0}^{s}e^{\tau - s}\,f_{8}(\tau\ d\tau\right)_{xx}ds, \vspace{0.2cm}  \\
(\rho_{2} + 4\,\gamma + 4\,\delta)\,v - b\,v_{xx} - 3\,k\,(\varphi_{x} + \psi -v) = (\rho_{2} + 4\,\gamma)\,f_{5} + \rho_{2}\,f_{6},
\end{array}
\right. 
\end{equation}
where 
\begin{equation*}
\tilde{g}_{i} = \int_{0}^{\infty} e^{-\,s}\,g_{i} (s)\ ds,\quad i=1,\,2.
\end{equation*}
We see that, if $\eqref{z5f5}$ admits a solution satisfying the required regularity in $\mathcal{D}(\mathcal{A})$, then \eqref{z1f1} implies that $\tilde{\varphi}$, $\tilde{\psi}$ and $\tilde{v}$ exist and satisfy the required regularity in $\mathcal{D}(\mathcal{A})$. On the other hand, \eqref{z7} implies that $\eta$ and $z$ exist and satisfy $\eta_{s} ,\,\eta\in L_{g_{1}}$ and $z_{s} ,\,z\in L_{g_{2}}$. Indeed, from $\eqref{z7*}$, we remark that it is enough to prove that $\eta\in L_{g_{1}}$ and $z\in L_{g_{2}}$. We have 
\begin{equation*}
s\mapsto (1 - e^{-\,s})(\varphi - f_{1})\in L_{g_{1}} \quad\hbox{and}\quad s\mapsto (1 - e^{-\,s})(\psi - f_{3})\in L_{g_{2}} 
\end{equation*}
because $\varphi,\,f_{1}\in H_{0}^{1} (0,\,1)$ and $\psi,\,f_{3}\in H_{*}^{1} (0,\,1)$. On the other hand, using the Fubini theorem and H\"older inequalities, we get
\begin{eqnarray*}
\lefteqn{\displaystyle\int_{0}^{1}\int_{0}^{\infty}g_1 (s)\left\vert\left( \displaystyle\int_{0}^{s}e^{\tau - s}\,f_{7}(\tau)\ d\tau \right)_x\right\vert^{2}ds\, dx }\nonumber \\ 
& \leq & \displaystyle\int_{0}^{\infty}e^{-\,2\,s}\,g_1 (s)\left(\displaystyle\int_{0}^{s}e^{\tau}\ d\tau\right)\displaystyle\int_{0}^{s}e^{\tau}\,\Vert f_{7x}(\tau)\Vert^{2}\, d\tau\, ds \\
\\
& \leq & \displaystyle\int_{0}^{\infty}e^{-\,s}\,(1 - e^{-\,s})\,g_1 (s)\displaystyle\int_{0}^{s}e^{\tau}\,\Vert f_{7x}(\tau)\Vert^{2}\,d\tau\,ds \\
\\
& \leq & \displaystyle\int_{0}^{\infty}e^{-\,s}\,g_1 (s)\displaystyle\int_{0}^{s} e^{\tau}\,\Vert f_{7x}(\tau)\Vert^{2}\, d\tau\, ds \\
\\
& \leq & \displaystyle\int_{0}^{\infty}e^{\tau}\,\Vert f_{7x}(\tau)\Vert^{2} \displaystyle\int_{\tau}^{+\infty}e^{-\,s}\,g_1 (s) \,ds\, d\tau \\
\\
& \leq & \displaystyle\int_{0}^{\infty}e^{\tau}\,g_1 (\tau)\,\Vert f_{7x}(\tau)\Vert^{2} \displaystyle\int_{\tau}^{+\infty}e^{-\,s}\,ds\,d\tau \\
\\
& \leq & \displaystyle\int_{0}^{\infty}g_{1}(\tau)\,\Vert f_{7x}(\tau)\Vert^{2}\, d\tau \\
\\
& \leq & \Vert f_{7}\Vert_{L_{g_1}}^{2} < \infty.
\end{eqnarray*} 
Then 
\begin{equation*}
s\mapsto \int_{0}^{s}e^{\tau - s}\,f_{7}(\tau)\ d\tau \in L_{g_{1}}.
\end{equation*}
Similarly, we get
\begin{equation*}
s\mapsto \int_{0}^{s}e^{\tau - s}\,f_{8}(\tau)\ d\tau \in L_{g_{2}}.
\end{equation*}
Therefore $\eta\in L_{g_{1}}$ and $z\in L_{g_{2}}$. Finally, to prove that $\eqref{z5f5}$ admits a solution satisfying the required regularity in $\mathcal{D}(\mathcal{A})$, we consider the variational formulation of $\eqref{z5f5}$ and using the Lax-Milgram theorem and classical elliptic regularity arguments. This proves that $\eqref{ZF}$ has a unique solution $U\in \mathcal{D}\left(\mathcal{A}\right)$. By the resolvent identity, we have $\lambda\,I - \mathcal{A}$ is surjective, for any $\lambda >0$ (see \cite{liu00}). Consequently, the Lumer-Phillips theorem implies that $\mathcal{A}$ is the infinitesimal generator of a C$_{0}$-semigroup of contractions on $\mathcal{H}$. 

\section{Stability}

In this section, we prove the stability of \eqref{222}. Let $U_{0} \in \mathcal{H}$, $U$ be the solution to \eqref{222} and $E$ be the energy of $U$ given by  
\begin{equation}\label{214}
E(t) = \frac{1}{2}\,\Vert U(t)\Vert_{\mathcal{H}}^{2} 
\end{equation}
\begin{equation*}
= \frac{1}{2}\left(3\,k\,\Vert \varphi_{x} +\psi - v\Vert^{2} + 3\,(b - g_{2}^{0})\Vert \psi_{x}\Vert^{2} - 3\,g_{1}^{0}\,\Vert \varphi_{x}\Vert^{2} + b\,\Vert v_{x}\Vert^{2} + 4\,\delta \Vert v \Vert^{2} \right)
\end{equation*}
\begin{equation*}
+ \frac{1}{2}\left(3\,\rho_{1}\,\Vert \varphi_{t}\Vert^{2} + 3\,\rho_{2}\, \Vert\psi_{t}\Vert^{2} + \rho_{2}\,\Vert v_{t} \Vert^{2} + 3\int_{0}^{\infty}\left(g_{1}(s)\Vert \eta_{x}\Vert^{2} + g_{2}(s)\Vert z_{x}\Vert^{2}\right)ds\right).
\end{equation*}
According to \eqref{222}, we have  
\begin{equation*}
E' (t)=\langle U_t (t)\, ,\, U(t)\rangle_{\mathcal{H}}=\langle \mathcal{A}U(t)\, ,\, U(t)\rangle_{\mathcal{H}} ,
\end{equation*}
where we use $'$ to denote the derivative with respect to $t$. So, using \eqref{dissp}, we find  
\begin{equation}
\label{213}
E'(t) = -\ 4\,\gamma\,\Vert v_{t}\Vert^{2} + \frac{1}{2}\int_{0}^{\infty}\left(g'_{1}(s)\Vert \eta_{x}\Vert^{2} + g'_{2}(s)\Vert z_{x}\Vert^{2}\right)ds \leq 0,
\end{equation} 
since $g_{1}$ and $g_{2}$ are nonincreasing.
\vskip0,1truecm
To state our stability result, we consider the following additional hypothesis on the relaxation functions $g_{1}$ and $g_{2}$: 
\vskip0,1truecm
{\bf (H2)} We assume that $g_{1}^{0} >0$, $g_{2}^{0} >0$ and there exist two positive constants $\alpha_{1}$ and $\alpha_{2}$ and an increasing strictly convex function 
$G: \mathbb{R}_{+}\rightarrow  \mathbb{R}_{+}$ of class $C^{1}(\mathbb{R}_{+})\cap C^{2}(0,\,\infty)$ satisfying 
\begin{equation}
\label{202}G(0) = G'(0) = 0\quad \mbox{and}\quad \lim_{t\rightarrow \infty}G'(t) = \infty
\end{equation}
such that, for any $i=1,\,2$,
\begin{equation}
\label{203+}g'_{i} (s) \leq -\, \alpha_{i}\,g_{i}(s),\quad s\in \mathbb{R}_{+}
\end{equation}
or
\begin{equation}
\label{203}\int_{0}^{\infty}\frac{g_{i}(s)}{G^{-1}\left(-\,g'_{i} (s)\right)}ds +\sup_{s\in\mathbb{R}_{+}}\frac{g_{i}(s)}{G^{-1}\left(-\,g'_{i}(s)\right)} < \infty.
\end{equation}
\vskip0,1truecm
\begin{theorem}
\label{theo:main} Assume {\bf (H1)} and {\bf (H2)} hold true. Let $U_{0} \in \mathcal{H}$ such that, for any $i=1,\,2$, 
\begin{equation}
\label{203*}\eqref{203+}\,\hbox{holds}\quad\hbox{or}\quad\sup_{t\in \mathbb{R}_{+}}\int_{t}^{\infty}\frac{g_i (s)}{G^{-1}\left(-\,g'_{i} (s)\right)}\Vert f_{0x}(\cdot ,\,s - t)\Vert^{2}\ ds < \infty,
\end{equation}
where $f_{0} =\varphi_{0}$ if $i=1$, and $f_0 =v_0 -u_0$ if $\,i=2$.
Then there exist two positive constants $c_{1}$ and $c_{2}$ such that the solution $U$ of \eqref{222} satisfies
\begin{equation}
\label{203**}E(t) \leq c_{2}\,G_{1}^{-\,1}(c_{1}\,t),\quad t\in \mathbb{R}_{+}, 
\end{equation}
where
\begin{equation}
\label{203***} G_{1}(s) = \int_{s}^{1}\frac{1}{G_{0}(\tau)}\ d\tau\quad\mbox{and}\quad G_{0}(s) =
\begin{cases}
s & \quad\hbox{if \eqref{203+} holds for any $i=1,\,2$}, \\
s\,G'(s) & \quad\hbox{otherwise}.
\end{cases}
\end{equation} 
\end{theorem}
\vskip0,1truecm
\begin{remark}
The hypothesis \eqref{203+} implies that $g_{i}$ converges exponentially to zero at infinity. So, if both $g_{1}$ and $g_{2}$ converge exponentially to zero at infinity, then \eqref{203**} leads to the exponential stability
\begin{equation}
\label{expon} E(t) \leq c_{2}\,e^{-\,c_{1}\,t},\quad t\in \mathbb{R}_{+} .
\end{equation} 
However, the hypothesis \eqref{203}, which was introduced by the second author in \cite{gues0}, allows $s\mapsto g_{i}(s)$ to have a decay rate at infinity arbitrarily closed to $\frac{1}{s}$. Indeed, for example, for $g_{i}(s) = d_{i}\,(1 + s)^{-\,q_{i}}$ with $d_{i} >0$ and $q_{i} >1$, hypothesis \eqref{203} is satisfied with $G(s)=s^{r}$, for all $r>\max\{\frac{q_{1} + 1}{q_{1} - 1},\frac{q_{2} + 1}{q_{2} - 1}\}$. And then \eqref{203**} implies that 
\begin{equation*}
E(t) \leq c_{2}\,(t + 1)^{\frac{-\,1}{r - 1}},\quad t\in \mathbb{R}_{+}.
\end{equation*}
For other examples, see \cite{gues0} and \cite{gues5}. In general, the decay rate of $E$ depends on the decay rate of $g_{i}$ which has the weaker decay rate. 
\end{remark}
\vskip0,1truecm
{\bf Proof of Theorem \ref{theo:main}.} In order to prove Theorem \ref{theo:main}, we will need to construct a Lyapunov functional equivalent to the energy $E$. To simplify the computations, we denote by $C$ and $C_{\lambda}$ some positive constants which can be different from line to line and depend contunuously on $E(0)$ and on some postive constant $\lambda$. In order to construct a Lyapunov functional $\mathcal{R}$ equivalent to the energy $E$, we define first some functionals $I_{1}$, $I_{2}$, $I_{3}$, $\mathcal{F}$ and $\mathcal{L}$ and prove some estimates on theirs derivatives. Let
\begin{equation}
\label{401} I_{1}(t) = -\ 3\,\rho_{1}\int_{0}^{1} \varphi_{t}\int_{0}^{\infty}g_{1}(s)\,\eta\ ds\ dx,\quad I_{2}(t) = -\ 3\,\rho_{2}\int_{0}^{1}\psi_{t}\int_{0}^{\infty}g_{2}(s)\,z\ ds\ dx,
\end{equation}
\begin{eqnarray}
\label{402}& & I_{3} (t) = \int_{0}^{1}\left(3\,\rho_{1}\,\varphi\,\varphi_{t} + 3\,\rho_{2}\,\psi\,\psi_{t} + \rho_{2}\,v\,v_{t}\right)dx,\quad \mathcal{F}(t)= I_{1}(t) + I_{2} (t) + I_{3}(t) \nonumber \\
& \hbox{and}& \quad \mathcal{L}(t) = E(t) + \varepsilon\,\mathcal{F}(t),
\end{eqnarray}
where $\varepsilon>0$ is a small parameter to be chosen later.
\vskip0,1truecm
\begin{lemma}
\label{lem:010}
For any $\delta_{0} >0$, there exists $C_{\delta_{0}}>0$ such that the functionals $I_{1}$ and $I_{2}$ satisfy	
\begin{eqnarray}
\label{4031} I'_{1} (t) \leq - 3\,\rho_{1}\,(g_{1}^{0} - \delta_{0})\Vert \varphi_{t}\Vert^{2} + \delta_{0}\left(\Vert \varphi_{x}\Vert^{2} + \Vert \varphi_{x} + \psi -v\Vert^{2}\right) + C_{\delta_{0}}\int_{0}^{\infty}g_{1}(s)\,\Vert\eta_{x}\Vert^{2}\  ds 
\end{eqnarray}
and 
\begin{eqnarray}
\label{4032} I'_{2} (t) \leq -\ 3\,\rho_{2}\,(g_{2}^{0} - \delta_{0})\,\Vert \psi_{t}\Vert^{2} + \delta_{0}\left(\Vert\psi_{x}\Vert^{2} + \Vert \varphi_{x} + \psi -v\Vert^{2}\right) + C_{\delta_{0}}\int_{0}^{\infty}g_{2}(s)\,\Vert z_{x}\Vert^{2}\ ds.
\end{eqnarray}
\end{lemma}
\vskip0,1truecm
{\bf Proof}. First, we note that 
\begin{equation*}
\begin{array}{lll}
\partial_{t}\displaystyle\int_{0}^{\infty}g_{1}(s)\,\eta\,ds & = & \partial_{t}\displaystyle\int_{-\infty}^{t}g_{1}(t - s)\,(\varphi (t) - \varphi (s))\ ds \vspace{0.2cm}\\
& = &\displaystyle\int_{-\infty}^{t}g^{\prime}_{1}(t - s)\,(\varphi (t) - \varphi (s))\ ds + \left(\displaystyle\int_{-\infty}^{t}g_{1}(t - s)\ ds\right)\varphi_{t}
\end{array}
\end{equation*}
that is
\begin{equation}
\label{eq2}
\partial_{t}\displaystyle\int_{0}^{\infty}g_{1}(s)\,\eta\ ds = \displaystyle\int_{0}^{\infty}g^{\prime}_{1}(s)\,\eta\ ds + g_{1}^{0}\,\varphi_{t}.
\end{equation}
Similarly, we have
\begin{equation}
\label{eq2+}
\partial_{t}\displaystyle\int_{0}^{\infty}g_{2}(s)\,z\ ds =\displaystyle\int_{0}^{\infty}g^{\prime}_{2}(s)\,z\ ds + g_{2}^{0}\,\psi_{t}.
\end{equation}
Second, using Young's and H\"{o}lder's inequalities, we get the following inequality: for all $\lambda >0$, there exists $C_{\lambda}>0$ such that, for any $v\in L^{2}(0,\,1)$ and ${\hat{\eta}} \in\{\eta,\,\eta_{x}\}$ in case $i=1$, and ${\hat{\eta}} \in\{z,\,z_{x}\}$ in case $i=2$, 
\begin{equation}
\label{in1}
\left\vert\displaystyle\int_{0}^{1}v\displaystyle\int_{0}^{\infty}g_{i}(s)\,{\hat {\eta}} \ ds\ dx\right\vert
\leq \lambda\,\Vert v\Vert^{2} + C_{\lambda}\displaystyle\int_{0}^{\infty}g_{i}(s)\, \Vert{\hat {\eta}}\Vert^{2}\ ds.
\end{equation}
Similarly,
\begin{equation}
\label{in2}
\left\vert\displaystyle\int_{0}^{1}v\displaystyle\int_{0}^{\infty }g_{i}^{\prime}(s)\,{\hat {\eta}}\ ds\ dx\right\vert
\leq \lambda\,\Vert v\Vert^{2} - C_{\lambda}\displaystyle\int_{0}^{\infty} g_{i}^{\prime}(s)\,\Vert{\hat {\eta}}\Vert^{2}\ ds.
\end{equation}
Now, direct computations, using the first equation in \eqref{109+}, integrating by parts and using the boundary conditions \eqref{114+} and \eqref{eq2}, yield 
\begin{equation*}
\begin{array}{lll}
I_{1}^{\prime }(t) & = & -\ 3\,\rho_{1}\,g_{1}^{0}\,\Vert\varphi_{t}\Vert^{2}  + 3\left\Vert\displaystyle\int_{0}^{\infty}g(s)\,\eta_{x}\ ds\right\Vert^{2}  - 3\,\rho_{1}\displaystyle\int_{0}^{1}\varphi_{t}\displaystyle\int_{0}^{\infty}g'_{1}(s)\,\eta\ ds\ dx \vspace{0.2cm} \\
& & \mbox{} + 3\,k\displaystyle\int_{0}^{1} (\varphi_{x} + \psi - v)\displaystyle\int_{0}^{\infty}g(s)\,\eta_{x}\ ds\ dx 
- 3\,g_{1}^{0}\displaystyle\int_{0}^{1}\varphi_{x}\displaystyle\int_{0}^{\infty}g(s)\,\eta_{x}\ ds\ dx.
\end{array}
\end{equation*}
Using \eqref{in1} and \eqref{in2} (with $i=1$) for the last three terms of this equality, Poincar\'{e}'s inequality for $\eta$, H\"older's inequality to estimate its second term, and \eqref{2.230} to estimate $-\,g_{1}^{\prime}$ by $\beta_{1} g_{1}$, we get \eqref{4031}. Similarly, using the second equation in \eqref{109+}, \eqref{eq2+}, \eqref{in1} and \eqref{in2} (with $i=2$), we find \eqref{4032}.
\vskip0,1truecm
\begin{lemma}
\label{lem:011}
For any $\delta_{0} >0$, there exists $C_{\delta_{0}}>0$ such that the functional $I_{3}$ satisfies	
\begin{equation}
\label{4033}I'_{3} (t) \leq -\ 3\,k\,\Vert \varphi_{x} + \psi - v\Vert^{2} + 3\,(g_{1}^{0} + \delta_{0})\Vert\varphi_{x}\Vert^{2} - 3\,(b - g_{2}^{0} - \delta_{0})\,\Vert\psi_{x}\Vert^{2} - (b - \delta_{0})\,\Vert v_{x}\Vert^{2} - 4\,\delta\,\Vert v\Vert^{2} 
\end{equation}
\begin{equation*}
+\ 3\,\rho_{1}\,\Vert \varphi_{t}\Vert^{2} + 3\,\rho_{2}\,\Vert \psi_{t}\Vert^{2} + C_{\delta_{0}}\,\Vert v_{t}\Vert^{2} + C_{\delta_{0}}\int_{0}^{\infty} (g_{1}(s)\,\Vert\eta_{x}\Vert^{2} + g_{2}(s)\,\Vert z_{x}\Vert^{2})\ ds.
\end{equation*}
\end{lemma}
\vskip0,1truecm
{\bf Proof}. Differentiating $I_{3}$, using the three equations in \eqref{109+}, integrating by parts and using the boundary conditions \eqref{114+}, we get
\begin{equation*}
I'_{3}(t) = -\ 3\,k\,\Vert \varphi_{x} + \psi - v\Vert^{2} + 3\,g_{1}^{0}\Vert\varphi_{x}\Vert^{2} - 3\,(b - g_{2}^{0})\,\Vert\psi_{x}\Vert^{2} - b\,\Vert v_{x}\Vert^{2} - 4\,\delta\,\Vert v\Vert^{2} + \rho_{2}\,\Vert v_{t}\Vert^{2}  
\end{equation*}
\begin{equation*}
+\ 3\,\rho_{1}\,\Vert \varphi_{t}\Vert^{2} + 3\,\rho_{2}\,\Vert \psi_{t}\Vert^{2} - 4\,\gamma\displaystyle\int_{0}^{1}v\,v_{t}\ dx - 3\displaystyle\int_{0}^{1} \varphi_{x}\int_{0}^{\infty}g_{1}(s)\,\eta_{x}\ ds\ dx - 3\displaystyle\int_{0}^{1} \psi_{x}\int_{0}^{\infty}g_{2}(s)\,z_{x}\ ds\ dx.
\end{equation*}
Applying \eqref{in1}  and Young's and Poincar\'e's inequalities to estimate the last three integrals of this equality, we obtain \eqref{4033}. 
\vskip0,1truecm
\begin{lemma}
\label{lem:02}
The functional $\mathcal{F}$ satisfies	
\begin{equation}
\label{403}\mathcal{F}'(t) \leq -\ E(t) + C\,\Vert v_{t}\Vert^{2} + C\int_{0}^{\infty} g_{1}(s)\,\Vert\eta_{x}\Vert^{2}\ ds + C\int_{0}^{\infty}g_{2}(s)\,\Vert z_x\Vert^{2} \ ds .
\end{equation}
Moreover, there exists a positive constant $\mu_{0}$ such that
\begin{equation}
\label{406}-\, \mu_{0}\,E \, \leq F \leq \mu_{0} \,E.
\end{equation}
\end{lemma}
\vskip0,1truecm
{\bf Proof}. By adding \eqref{4031}, \eqref{4032} and \eqref{4033}, and using the definition of $E$, we have
\begin{equation*}
\mathcal{F}'(t) \leq -\ 2\,E(t) + \delta_{0}\,(2\,\Vert \varphi_{x} +\psi - v\,\Vert^{2} + 4\,\Vert \varphi_{x} \Vert^{2} + 4\,\Vert \psi_{x} \Vert^{2} + \Vert v_{x} \Vert^{2} + 3\,\rho_{1}\,\Vert \varphi_{t} \Vert^{2} + 3\,\rho_{2}\,\Vert \psi_{t} \Vert^{2})
\end{equation*}
\begin{equation*}
+\ C_{\delta_{0}}\,\Vert v_{t}\Vert^{2} + C_{\delta_{0}}\int_{0}^{\infty} (g_{1}(s)\,\Vert\eta_{x}\Vert^{2} + g_{2}(s)\,\Vert z_{x}\Vert^{2} )\ ds.
\end{equation*}
On the other hand, the Poincar\'e's inequality applied to $\psi$ and $v$, and the condition \eqref{2.16} imply that
\begin{equation*}
\delta_{0}\,(2\,\Vert \varphi_{x} + \psi - v\Vert^{2} + 4\,\Vert \varphi_{x} \Vert^{2} + 4\,\Vert \psi_{x} \Vert^{2} + \Vert v_{x} \Vert^{2} + 3\,\rho_{1}\,\Vert \varphi_{t} \Vert^{2} + 3\,\rho_{2}\Vert \psi_{t} \Vert^{2}) \leq \delta_{0}\,C\,E(t).
\end{equation*}
Then combining the last two inequalities and choosing $\delta_0 >0$ small enough such that $\delta_{0}\,C\leq 1$, we find \eqref{403}. 
\vskip0,1truecm
From the Young's and Poincar\'e's inequalities, and the condition \eqref{2.16} lead to, for some $\mu_{0} >0$,
\begin{equation*}
\vert F\vert\leq \vert I_{1}\vert + \vert I_{2}\vert + \vert I_{3}\vert\leq \mu_{0} \,E,
\end{equation*}
\vskip0,1truecm
\begin{lemma}
\label{lem:02+}
There exits a positive constant $\varepsilon$ such that the functional $L$ satisfies	
\begin{equation}
\label{304}\mathcal{L}' (t)  \leq  -\ \varepsilon\,E(t) + C\int_{0}^{\infty}g_{1}(s) \,\Vert \eta_{x}\Vert^{2}\ ds + C\int_{0}^{\infty}g_{2}(s)\,\Vert z_{x}\Vert^{2}\ ds.
\end{equation}
and there exist positive constants $\mu_{1}$ and $\mu_{2}$ such that
\begin{equation}
\label{410} \mu_{1}\,E \, \leq L \leq \mu_{2} \,E.
\end{equation}
\end{lemma}
\vskip0,1truecm
{\bf Proof}. We differentiate $L$ from \eqref{402} with respect to time and use \eqref{403}, together with the dissipation of energy \eqref{213} and the nonincreasingness of $g_{1}$ and $g_{2}$, we obtain
\begin{equation*}
\mathcal{L}' (t) \leq -\ \varepsilon\,E(t) + (\varepsilon\,C - 4\,\gamma)\,\Vert v_{t}\Vert^{2} + \varepsilon\,C\int_{0}^{\infty}g_{1}(s)\,\Vert \eta_{x}\Vert^{2}\ ds +\varepsilon\,C\int_{0}^{\infty}g_{2}(s)\,\Vert z_{x}\Vert^{2}\ ds.
\end{equation*}
Then, for 
\begin{equation}
\label{3040} 0<\varepsilon <\frac{4\,\gamma}{C},
\end{equation}
\eqref{304} holds. On the other hand, thanks to \eqref{406}, we get
\begin{equation*}
(1 - \varepsilon\,\mu_{0})\,E \leq L \leq (1 + \varepsilon\, \mu_{0} )\,E.
\end{equation*}
So, for 
\begin{equation}
\label{3040+} 0<\varepsilon <\frac{1}{\mu_{0}},
\end{equation}
we find \eqref{410} with $\mu_{1} =1 - \varepsilon\,\mu_{0}$ and $\mu_{2} =1 + \varepsilon\,\mu_{0}$. Finally, choosing
\begin{equation*}
0<\varepsilon <\min\left\{\frac{4\,\gamma}{C},\ \frac{1}{\mu_{0}}\right\} ,
\end{equation*}  
we get \eqref{304} and \eqref{410}. 
\vskip0,1truecm
To estimate the last two terms of \eqref{304}, we adapt to our system a lemma introduced by the first author in \cite{gues0} and improved in \cite{gues9}.  
\vskip0,1truecm
\begin{lemma}
\label{lem:020} 
There exist positive constants $d_{1}$ and $d_{2}$ such that, for any $\varepsilon_{0}>0$, the following two inequalities hold:
\begin{equation}
\frac{G_{0}\,(\varepsilon_{0}\,E(t))}{\varepsilon_{0}\,E(t)}\displaystyle\int_{0}^{\infty}g_{1}(s)\,\Vert\eta_{x}\Vert^{2}\ ds \leq -\ d_{1}\,E^{\prime}(t) +  d_{1}\,G_{0}\,(\varepsilon_{0}\,E(t))\label{3.23}
\end{equation}
and
\begin{equation}
\frac{G_{0}\,(\varepsilon_{0}\,E(t))}{\varepsilon_{0}\,E(t)}\displaystyle\int_{0}^{\infty}g_{2}(s)\,\Vert z_{x}\Vert^{2}\ ds \leq -\ d_{2}\,E^{\prime}(t) +  d_{2}\,G_{0}\,(\varepsilon_{0}\,E(t)).\label{3.23+}
\end{equation}
\end{lemma}
\vskip0,1truecm
{\bf Proof}. 
If \eqref{203+} holds, then we have from \eqref{213} 
\begin{equation}
\label{408}\int_{0}^{\infty}g_{i}(s)\,\Vert f_{x}\Vert^{2}\ ds \leq -\ \frac{1}{\alpha_{i}}\int_{0}^{\infty}g'_{i}(t)\,\Vert f_{x}\Vert^{2}\ ds \leq -\, \frac{2}{\alpha_{i}}\,E' (t),
\end{equation}
where $f=\eta$ in case $i=1$, and $f=z$ in case $i=2$. So \eqref{3.23} and \eqref{3.23+} hold with $d_{i} = \frac{2}{\alpha_{i}}$ and $G_{0}(s)=s$.
\vskip0,1truecm
When \eqref{203} is satisfied, we note first that, if $E (t_{0})=0$, for some $t_{0} \in \mathbb{R}_+$, then $E(t)=0$, for all $t\geq t_{0}$, since $E$ is nonnegative and nonincreasing, and consequently, \eqref{203**} is satisfied, since $E$ is bounded. Thus, without loss of generality, we can assume that $E>0$ on $\mathbb{R}_+$.
\vskip0,1truecm
Because $E$ is nonincreasing, we have
\begin{equation*}
\Vert\eta_{x}\Vert^{2} \leq 2\left(\Vert\varphi_{x}(\cdot,\,t)\Vert^{2} + \Vert\varphi_{x} (\cdot,\,t - s)\Vert^{2} \right)  
\leq C\,E(0) + 2\,\Vert\varphi_{x}(\cdot,\,t - s)\Vert^{2}  
\end{equation*}
\begin{equation*} 
\leq \left\{
\begin{array}{ll}
C\,E(0) & \quad\hbox{if} \quad 0\leq s\leq t, \\
C\,E(0) + 2\Vert\varphi_{0x} (\cdot,\,s - t)\Vert^{2}&\quad\hbox{if}\quad s>t\geq 0
\end{array} :=M_{1}(t,\,s), 
\right.
\end{equation*}
so we conclude that 
\begin{equation}
\Vert\eta_{x}\Vert^{2} \leq M_{1}(t,\,s),\quad t,\,s\in \mathbb{R}_{+}. \label{3.26} 
\end{equation} 
Similarly, we get
\begin{equation}\label{3.260}
\Vert z_{x}\Vert^{2} \leq \left\{
\begin{array}{ll}
C\,E(0) & \ \hbox{if} \quad 0\leq s\leq t, \\
C\,E(0) + 2\,\Vert v_{0x} (\cdot,\,s - t) - u_{0x}(\cdot,\,s - t)\Vert^{2} & \ \hbox{if}\,\, s>t\geq 0
\end{array} :=M_{2}(t,\,s).
\right.
\end{equation} 
Let $\tau_{1}(t,\,s),\ \tau_{2}(t,\,s)>0$ (which will be fixed later on), $\varepsilon_{0} >0$ and $K(s)={\frac{s}{{G^{-\,1}(s)}}}$, for $s>0$, and $K(0)=0$, since {\bf (H2)} implies that  
\begin{equation*}
\lim_{s\to 0^{+}}\dfrac{s}{G^{-1}(s)}=\lim_{\tau\to 0^{+}}\dfrac{G(\tau)}{\tau}=G'(0)=0.
\end{equation*}
The function $K$ is nondecreasing. Indeed, the fact that $G^{-1}$ is concave and $G^{-1} (0)=0$ implies that, for any 
$0\leq s_{1} < s_{2}$,
\begin{equation*}
K(s_{1}) = \dfrac{s_{1}}{G^{-\,1}\left(\frac{s_{1}}{s_{2}}s_{2} + \left(1 - \dfrac{s_{1}}{s_{2}}\right)0\right)} 
\leq \dfrac{s_{1}}{\dfrac{s_{1}}{s_{2}}\,G^{-1}(s_{2}) + \left(1 - \dfrac{s_{1}}{s_{2}}\right)G^{-1} (0)} = \dfrac{s_{2}}{G^{-1}(s_{2})} = K(s_{2} ).
\end{equation*} 
Then, using \eqref{3.26} and \eqref{3.260},
\begin{equation}
K\left(-\,\tau_{2} (t,\,s)\,g_{1}^{\prime}(s)\,\Vert \eta_{x}\Vert^{2}\right) \leq K\left(-\,M_{1}(t,\,s)\,\tau_{2}(t,\,s)\,g_{1}^{\prime}(s)\right),\quad t,\,s\in \mathbb{R}_{+} \label{3.2601} 
\end{equation}
and
\begin{equation}
K\left(-\,\tau_{2} (t,\,s)\,g_{2}^{\prime}(s)\Vert z_{x}\Vert^{2}\right) \leq K\left(-\,M_{2} (t,\,s) \,\tau_{2}(t,\,s)\,g_{2}^{\prime}(s)\right),\quad t,\,s\in \mathbb{R}_{+} . \label{3.2602} 
\end{equation}
Using \eqref{3.2601}, we arrive at 
\begin{eqnarray*}
\int_{0}^{\infty}g_{1}(s)\,\Vert \eta_{x}\Vert^{2}\ ds = \dfrac{1}{G'(\varepsilon_{0}\,E(t))}\int_{0}^{\infty} \dfrac{1}{\tau_{1}(t,\,s)}G^{-\,1}\left(-\,\tau_{2}(t,\,s)\,g'_{1}(s)\,\Vert \eta_{x}\Vert^{2}\right) 
\end{eqnarray*} 
\begin{equation*}
\times\dfrac{\tau_{1}(t,\,s)\,G'(\varepsilon_{0}\,E(t))\,g_1 (s)}{-\,\tau_{2}(t,\,s)\,g'_1(s)}\ K\left(-\,\tau_{2}(t,\,s)\,g'_{1}(s)\Vert \eta_{x}\Vert^{2}\right)ds
\end{equation*} 
\begin{equation*}   
\leq \dfrac{1}{G'(\varepsilon_{0}\,E(t))}\int_{0}^{\infty}\dfrac{1}{\tau_{1}(t,\,s)}\,G^{-\,1}\left(-\,\tau_{2}(t,\,s)\,g'_{1}(s)\,\Vert \eta_{x}\Vert^{2}\right)\dfrac{\tau_{1}(t,\,s)G'(\varepsilon_{0}\,E(t))\,g_{1}(s)}{-\,\tau_{2}(t,\,s)\,g'_1 (s)}\ K\left(-\,M_{1}(t,\,s)\,\tau_{2}(t,\,s)\,g'_1 (s)\right)ds
\end{equation*} 
\begin{equation*}
\leq\dfrac{1}{G' (\varepsilon_{0}\,E(t))}\int_{0}^{\infty}\dfrac{1}{\tau_{1}(t,\,s)}\,G^{-1}\left(-\,\tau_{2}(t,\,s)\,g'_{1}(s)\,\Vert \eta_{x}\Vert^{2}\right) \dfrac{M_{1}(t,\,s)\,\tau_{1}(t,\,s)\,G'(\varepsilon_{0} \,E(t))\,g_{1}(s)}{G^{-\,1}\,(-\,M_{1}(t,\,s)\,\tau_{2}(t,\,s)\,g'_{1}(s))}\ ds.
\end{equation*}
Let $G^{*} (s)=\sup_{\tau\in \mathbb{R}_{+}} \{s\,\tau - G(\tau)\}$, for $s\in \mathbb{R}_{+}$, denote the dual function of $G$. Thanks to {\bf (H2)}, we see that
\begin{equation*}
G^{*}(s) = s\,(G')^{-1} (s) - G ((G')^{-1} (s)),\quad s\in \mathbb{R}_{+} .
\end{equation*}
Using Young's inequality: $s_{1}\,s_{2}\leq G(s_{1}) + G^{*}(s_{2})$, for 
\begin{equation*}
s_{1} = G^{-1}\left(-\,\tau_{2}(t,\,s)\,g_1^{\prime}(s)\Vert \eta_x\Vert^{2}\right)\quad\mbox{and}\quad s_{2} ={\dfrac{{M_{1}(t,\,s)\,\tau_{1}(t,\,s)\,G^{\prime}(\varepsilon_{0}\,E(t))\,g_1 (s)}}{{G^{-\,1}(-\,M_{1}(t,\,s)\,\tau_{2}(t,\,s) g_1^{\prime }(s))}}}, 
\end{equation*}
we get 
\begin{equation*}
\int_{0}^{\infty}g_{1}(s)\Vert \eta_{x}\Vert^{2}\ ds \leq  \dfrac{1}{G^{\prime}(\varepsilon_{0}\,E(t))}\int_{0}^{\infty}\dfrac{{-\,\tau_{2}(t,\,s)}}{{\tau_{1}(t,\,s)}}g_{1}^{\prime}(s)\Vert \eta_{x}\Vert^{2}\ ds 
\end{equation*} 
\begin{equation*}
+\ \dfrac{1}{G^{\prime}(\varepsilon_{0}\,E(t))}\int_{0}^{\infty}\dfrac{1}{\tau_{1}(t,\,s)}\ G^{*}\left({\frac{{M_{1}(t,\,s)\,\tau_{1}(t,\,s)\,G^{\prime}(\varepsilon_{0}\,E(t))\,g_1 (s)}}{{G^{-1}(-\,M_{1}(t,\,s)\,\tau_{2}(t,\,s)\,g^{\prime}_1 (s))}}}\right)ds. 
\end{equation*}
Using the fact that $G^{*}(s)\leq s\,(G')^{-\,1} (s)$, we get 
\begin{equation*}
\int_{0}^{\infty}g_1 (s)\Vert \eta_x\Vert^{2}\, ds  \leq  \dfrac{-\,1}{G^{\prime}(\varepsilon_{0}\,E(t))}\int_{0}^{\infty}\dfrac{{\tau_{2}(t,\,s)}}{{\tau_{1}(t,\,s)}}\ g_{1}^{\prime}(s)\Vert \eta_{x}\Vert^{2}\ ds 
\end{equation*}
\begin{equation*}
\mbox{} + \int_{0}^{\infty}{\frac{{M_{1}(t,\,s)\,g_{1}(s)}}{{G^{-\,1}(-\,M_{1}(t,\,s)\,\tau_{2}(t,\,s)\,g_{1}^{\prime}(s))}}}\,(G')^{-\,1}\left({\frac{{M_1 (t,\,s)\,\tau_{1}(t,\,s)\,G^{\prime}(\varepsilon_{0}\,E(t))\,g_1 (s)}}{{G^{-\,1}(-\,M_1 (t,\,s)\,\tau_{2}(t,\,s)\,g_1^{\prime}(s))}}}\right)ds. 
\end{equation*}
Then, using the fact that $(G')^{-1}$ is nondecreasing and choosing $\tau_{2}(t,\,s) = {\frac{{1}}{{M_{1}(t,\,s)}}}$, we get 
\begin{equation*}
\int_{0}^{\infty}g_{1}(s)\,\Vert \eta_{x}\Vert^{2}\ ds \leq \dfrac{-\,1}{G^{\prime}(\varepsilon_{0}\,E(t))}\int_{0}^{\infty} \dfrac{1}{M_{1}(t,\,s)\,\tau_{1}(t,\,s)}\,g_{1}^{\prime}(s)\,\Vert \eta_{x}\Vert^{2}\ ds  
\end{equation*} 
\begin{equation*}
+ \int_{0}^{\infty}{\dfrac{{M_{1}(t,\,s)\,g_{1}(s)}}{{G^{-\,1}(-\,g_{1}^{\prime}(s))}}}\ (G')^{-1}\left(m_{1}\,M_{1}(t,\,s)\,\tau_{1}(t,\,s)\,G^{\prime}(\varepsilon_{0}\,E(t))\right)ds,  
\end{equation*}
where $m_{1} = \sup_{s\in \mathbb{R}_{+}}{\frac{{g_{1}(s)}}{{G^{-1}(-\,g^{\prime}_{1}(s))}}} <\infty$ ($m_{1}$ exists according to \eqref{203}). Due to \eqref{203} and the restriction on $\varphi_{0}$ in \eqref{203*}, we have
\begin{equation*}
\sup_{t\in \mathbb{R}_{+}}\displaystyle\int_0^{\infty}{\frac{{M_{1}(t,\,s)\,g_{1}(s)}}{{G^{-\,1}(-\,g_{1}^{\prime}(s))}}}\ ds=:m_{2} <\infty .
\end{equation*} 
Therefore, choosing $\tau_{1}(t,\,s)= {\frac{{1}}{{m_{1}\,M_{1}(t,\,s)}}}$ and using \eqref{213}, we obtain 
\begin{equation*}
\int_{0}^{\infty}g_{1}(s)\Vert \eta_{x}\Vert^{2}\ ds \leq \dfrac{{-\,m_{1}}}{{G^{\prime}(\varepsilon_{0}\,E(t))}}\int_{0}^{\infty}g_{1}^{\prime}(s)\,\Vert \eta_{x}\Vert^{2}\ ds + \varepsilon_{0}\,E(t)\int_{0}^{\infty}{\frac{{M_{1}(t,\,s)\,g_{1}(s)}}{{G^{-\,1}(-\,g_{1}^{\prime}(s))}}}\ ds  
\end{equation*} 
\begin{equation*}
\leq \dfrac{{-\,2\,m_{1}}}{{G^{\prime}(\varepsilon_{0}\,E(t))}}\,E^{\prime}(t) + \varepsilon_{0}\,m_{2}\,E(t),
\end{equation*}
which gives \eqref{3.23} with $d_{1} = \max \{2\,m_{1},\,m_{2}\}$ and $G_{0}(s) = s\,G' (s)$. Repeating the same arguments, we get
\eqref{3.23+} with $G_{0}(s) = s\,G' (s)$, $d_{2} = \max \{2\,m_{1},\,m_{2}\}$,  
\begin{eqnarray*}
m_{1} = \sup_{s\in \mathbb{R}_{+}}{\frac{{g_{2}(s)}}{{G^{-1}(-\,g_{2}^{\prime}(s))}}} ,\quad m_{2} =\sup_{t\in \mathbb{R}_{+}}\displaystyle\int_{0}^{\infty}{\frac{{M_{2}(t,\,s)\,g_{2} (s)}}{{G^{-\,1}(-\,g_{2}^{\prime}(s))}}}\ ds, \\
\tau_{1}(t,\,s)= {\frac{{1}}{{m_{1}\,M_{2}(t,\,s)}}}\quad\hbox{and}\quad \tau_{2} (t,\,s) =  {\frac{{1}}{{M_{2}(t,\,s)}}}.
\end{eqnarray*}
\vskip0,1truecm
We are now ready to prove the main stability result \eqref{203**}. Multiplying \eqref{304} by $\frac{G_{0}(\varepsilon_{0}\,E)}{\varepsilon_{0}\,E}$ and combining \eqref{3.23} and \eqref{3.23+}, 
we find, for $c_{0} =d_{1} + d_{2}$,
\begin{equation}
\label{304**}
\frac{G_{0}(\varepsilon_{0}\,E(t))}{\varepsilon_{0}\,E(t)}\,\mathcal{L}' (t)  + c_{0}\,C\,E'(t) \leq -\,\left(\frac{\varepsilon}{\varepsilon_{0}} - c_{0}\,C \right)G_{0}(\varepsilon_{0}\,E(t)).
\end{equation}
We define our Lyapunov functional $\mathcal{R}$ by  
\begin{equation}
\label{304*}
\mathcal{R} = \tau_{0} \left(\frac{G_{0}(\varepsilon_{0}\,E)}{\varepsilon_{0}\,E}\mathcal{L} + c_{0}\,C\,E\right),
\end{equation}
where $\tau_{0}$ is a positive constant that will be choosen later. Because $E$ is nonincreasing and $G$ is convex, then $\frac{G_{0}(\varepsilon_{0}\,E)}{\varepsilon_{0}\,E}$ is nonincreasing, and therefore, \eqref{304**} and \eqref{304*} lead to 
\begin{equation}\label{3041}
\mathcal{R}' (t)  \leq -\,\tau_{0}\,\left(\frac{\varepsilon}{\varepsilon_{0} } - c_{0}\,C \right)G_{0}(\varepsilon_{0}\,E(t)).
\end{equation}
Moreover, recalling that $\frac{G_{0}(\varepsilon_{0}\,E)}{\varepsilon_{0}\,E}$ is nonincreasing and using \eqref{410}, we obtain 
\begin{equation}
\label{4101} \tau_{0}\,c_{0}\,C\,E \leq \mathcal{R} \leq \tau_{0}\left[(1 + \varepsilon\,\mu_{0})\,\frac{G_{0}(\varepsilon_{0}\,E(0))}{\varepsilon_{0}\,E(0)} + c_{0}\,C\right]E.
\end{equation}
By choosing $0<\varepsilon_{0} <\frac{\varepsilon}{c_{0}\,C}$, we deduce from \eqref{3041} and \eqref{4101} that $\mathcal{R}$ is equivalent to $E$ and satisfies  
\begin{equation}
\label{3042}
\mathcal{R}' (t)  \leq -\,\tau_{0}\,C\,G_{0}(\varepsilon_{0}\,E(t)).
\end{equation}
Thus, for $\tau_{0} >0$ such that
\begin{equation*}
\mathcal{R} \leq \varepsilon_{0}\,E \quad\hbox{and}\quad \mathcal{R}(0) \leq 1,
\end{equation*} 
we get, for $c_{1} = \tau_{0}\,C$,
\begin{equation}
\label{3.30}
\mathcal{R}^{\prime} (t) \leq -\ c_{1}\,G_{0}( \mathcal{R}(t)).
\end{equation}
Then \eqref{3.30} implies that $(G_{1}(\mathcal{R}))' \geq c_{1}$, where $G_{1}$ is defined in \eqref{203***}. So, a direct integrating gives 
\begin{equation}
\label{3.300}
G_{1}(\mathcal{R} (t)) \geq c_{1}\,t + G_{1}(\mathcal{R} (0)).
\end{equation}
Because $\mathcal{R} (0)\leq 1$ and $G_{1}$ is decreasing, we obtain $G_{1} (\mathcal{R} (t)) \geq c_{1}\,t$ which implies that $\mathcal{R} (t)\leq G_{1}^{-1} (c_{1}\,t)$. Finally, from the equivalence of $\mathcal{R}$ and $E$, the result \eqref{203**} follows and the proof of Theorem \ref{theo:main} is complete.

\section{Numerical analysis}

\setcounter{equation}{0}
In this section, we present some numerical results illustrating the asymptotic behavior of the energy for the exponential decay. For this, we use Finite Difference (of
second order in space and time). Furthermore, the method of $\beta-$Newmark is a second order method preserving the discrete energy always when the discrete system of equations of motion is symmetric (i.e. matrices associated to the system should be symmetric).

\subsection{Finite difference method.}
We consider $J$ an integer non-negative and $h = L/(J + 1)$ an spatial subdivision of the interval $(0, L)$ given by
$0=x_0 <x_1 <\ldots<x_J <x_{J+1} =L$, 
with $x_j = jh$ each node of the mesh. 
We use $\varphi_j(t)$,  $\psi_j(t)$, $v_j(t)$, for all $j = 1,2,\ldots,J$ and $t > 0$ to denote the approximate values of $\varphi(jh, t)$ and $\psi(jh, t)$, respectively.
In addition, we denote the discrete operator $\Delta_{h}\vartheta_j=\frac{\vartheta_{j+1}-2\vartheta_j+\vartheta_{j-1}}{h^2}$,
$\delta^-_h\vartheta_j=\frac{\vartheta_{j}-\vartheta_{j-1}}{h}$
and the $\theta$-scheme
$Q_\theta\vartheta_j={\theta\vartheta_{j+1}+(1-2\theta)\vartheta_j+\theta\vartheta_{j-1}}$,
with $\theta=1/4$.
We assume the following finite difference scheme applied to system
\begin{align}\label{109bb}
\begin{cases}
3\left(\rho_{1}\,\varphi''_{j} - k\,(\Delta_h\varphi_j +\delta^-_h\psi_j -\delta^-_h v_j) + \displaystyle\int_0^{\infty}\,g_1 (s) \Delta_h\varphi_{j} (t-s) \,ds\right) = 0, \quad & j=1,\ldots,J,\\
3\left(\rho_{2}\,\psi''_{j} - (b-g_2^0 )\,\Delta_h\psi_{j} + k\,(\delta^-_h\varphi_{j} +
Q_{1/4}\psi_j -Q_{1/4}v_j)  + \displaystyle\int_0^{\infty}\,g_2 (s)\Delta_h\psi_{j} \,ds\right) =0, &j=1,\ldots,J,\\
\rho_{2}\,v''_{j} - b\,\Delta_hv_{j} - 3\,k\,(\delta^-_h\varphi_{j} +Q_{1/4}\psi_j -Q_{1/4}v_j) + 4\,\Theta_{1/4}\delta\,v_j +4\gamma v'_j= 0,&
j=1,\ldots,J,
\\
\varphi_0=\varphi_J=\psi_0=v_0=0,\quad  \psi_J=\psi_{J+1},\quad v_J=v_{J+1}.&
\end{cases}
\end{align}

\subsection{Equation of motion and time discretization}

The system \eqref{109bb} can be rewritten as
\begin{eqnarray}
\mathbf{M}\left[\begin{array}{c}
\ddot{\varphi}_h\\ \ddot{\psi}_h \\ \ddot{v}_h
\end{array}\right]
\,+\,
\mathbf{C}\left[\begin{array}{c}
\dot{\varphi}_h\\ \dot{\psi}_h \\ \dot{v}_h
\end{array}\right]
\,+\,
\mathbf{K}\left[\begin{array}{c}
{\varphi}_h\\ {\psi}_h \\ v_h
\end{array}\right]
\,+\,
{\mathbf{G}}\left(\left[\begin{array}{c}
{\varphi}_h\\ {\psi}_h \\ v_h
\end{array}\right]\right)
&=&
\mathbf{0},
\label{401}
\end{eqnarray}
where 
$\phi_h=(\phi_1,\ldots,\phi_J)$, $\psi_h=(\psi_1,\ldots,\psi_J)$, 
$v_h=(v_1,\ldots,v_J)\in \R^J $.
The mass, damping, and stiffness matrices  
($\mathbf{M}$, $\mathbf{C}$, $\mathbf{K} \in \mathcal{M}_{3J}(\R)$),
are given by 
\\
$\mathbf{M}=\begin{bmatrix}
3\rho_1 \mathbf{I} & \mathbf{0} & \mathbf{0} \\ \mathbf{0} & 3\rho_2  \mathbf{I} & \mathbf{0} \\ \mathbf{0} & \mathbf{0} & \rho_2  \mathbf{I}
\end{bmatrix}$, 
$\mathbf{C}=\begin{bmatrix}
\mathbf{0} & \mathbf{0} & \mathbf{0} \\ \mathbf{0} & \mathbf{0} & \mathbf{0} \\ \mathbf{0} & \mathbf{0} & 4 \gamma  \mathbf{I}
\end{bmatrix}$,
where $\mathbf{I}$ is the identity and  $\mathbf{0}$ is the null matrix of $\mathcal{M}_{J}(\R)$.
$\mathbf{K}=\begin{bmatrix}
- 3 k\mathbf{D}^2_0 & - 3 k\mathbf{D}^- & 3 k\mathbf{D}^- \\ 
3 k\mathbf{D}^- & - 3 b\mathbf{D}^+\mathbf{D}^- + 3 k\mathbf{Q} & -3 k\mathbf{Q} \\ 
-3k\mathbf{D}^- & -3k\mathbf{Q} & -b\mathbf{D}^+\mathbf{D}^- + (3k+4\delta)\mathbf{Q}
\end{bmatrix}$, 
where, 
$\mathbf{D}^2_0=\frac{1}{h^2}\begin{pmatrix}
-2 & 1 & \ldots & \ldots \\
1 & -2 & 1 & \ldots \\
\ldots & \ldots & \ldots & \ldots \\
\ldots & \ldots & 1 & -2
\end{pmatrix}$,
$\mathbf{D}^-=\frac{1}{h}\begin{pmatrix}
1 & 0 & \ldots & \ldots \\
-1 & 1 & 0 & \ldots \\
\ldots & \ldots & \ldots & \ldots \\
\ldots & \ldots & -1 & 1
\end{pmatrix}$, and $\mathbf{D}^+=- \left(\mathbf{D}^-\right)^t$.
We note that 
\begin{equation}
\mathbf{D}^2_0=\mathbf{D}^+\mathbf{D}^- - \frac{1}{h^2}\mathbf{J}^{JJ}
= \mathbf{D}^-\mathbf{D}^+ - \frac{1}{h^2}\mathbf{J}^{11},
\label{eq43b}
\end{equation}
where $\mathbf{J}^{ij}$ is the single-entry matrix of one, for $i$-row, $j$-column.
Additionally,\\
$\mathbf{Q}=\frac{1}{4}\begin{pmatrix}
2 & 1 & \ldots & \ldots \\
1 & 2 & 1 & \ldots \\
\ldots & \ldots & \ldots & \ldots \\
\ldots & \ldots & 1 & 1
\end{pmatrix}=
\mathbf{P}^+\mathbf{P}^-$, where
$\mathbf{P}^-=\frac{1}{2}\begin{pmatrix}
1 & 0 & \ldots & \ldots \\
1 & 1 & 0 & \ldots \\
\ldots & \ldots & \ldots & \ldots \\
\ldots & \ldots & 1 & 1
\end{pmatrix}$ and $\mathbf{P}^+=(\mathbf{P}^-)^t$.
The memory terms are given by,
$
{\mathbf{G}}\left(\left[\begin{array}{c}
{\varphi}_h\\ {\psi}_h \\ v_h
\end{array}\right]\right)=
\begin{pmatrix}
{\mathbf{G}}_1({\varphi}_h) \\ {\mathbf{G}}_2({\psi}_h) \\ 0
\end{pmatrix},
$
where
${\mathbf{G}}_1({\varphi}_h)=
\int_0^{\infty}\,g_1 (s) \mathbf{D}_0^2\varphi_{j} (t-s) \,ds $
and
${\mathbf{G}}_2({\psi}_h)=
\int_0^{\infty}\,g_2 (s) \mathbf{D}^+\mathbf{D}^-\psi_{j} (t-s) \,ds $.
\vskip0,1truecm
The Newmark algorithm \cite{Newmark}
is based on a set of two relations expressing the forward displacement
$[\varphi_h^{n+1}, \psi_h^{n+1}, v_h^{n+1}]^\top$
and velocity
$[\dot{\varphi}_h^{n+1}, \dot{\psi}_h^{n+1},  \dot{v}_h^{n+1}]^\top$.
The method  consists in updating
the displacement, velocity and acceleration vectors
from current time $t^n=n\delta t$ to the time $t^{n+1} = (n+1)\delta t$,
\begin{eqnarray}
\dot{\varphi}_{h}^{n+1} &=& \dot{\varphi}_{h}^{n} + (1-\varsigma)\delta t\,\ddot{\varphi}_{h}^{n} + \varsigma\delta t\,\ddot{\varphi}_{h}^{n+1}\label{402}\\
{\varphi}_{h}^{n+1} &=& {\varphi}_{h}^{n} + \delta t \dot{\varphi}_{h}^{n}
+ \left(\frac{1}{2}-\beta\right)\delta t^ 2\,\ddot{\varphi}_{h}^{n} + \beta\delta t^2\,\ddot{\varphi}_{h}^{n+1}\label{403}\\
\dot{\psi}_{h}^{n+1} &=& \dot{\psi}_{h}^{n} + (1-\varsigma)\delta t\,\ddot{\psi}_{h}^{n} + \varsigma\delta t\,\ddot{\psi}_{h}^{n+1}\label{404}\\
{\psi}_{h}^{n+1} &=& {\psi}_{h}^{n} +  \delta t \dot{\psi}_{h}^{n}
+ \left(\frac{1}{2}-\beta\right)\delta t^ 2\,\ddot{\psi}_{h}^{n} + \beta\delta t^2\,\ddot{\psi}_{h}^{n+1},\label{405}\\
\dot{v}_{h}^{n+1} &=& \dot{v}_{h}^{n} + (1-\varsigma)\delta t\,\ddot{v}_{h}^{n} + \varsigma\delta t\,\ddot{v}_{h}^{n+1}\label{404b}\\
v_{h}^{n+1} &=& v_{h}^{n} +  \delta t \dot{v}_{h}^{n}
+ \left(\frac{1}{2}-\beta\right)\delta t^ 2\,\ddot{v}_{h}^{n} + \beta\delta t^2\,\ddot{v}_{h}^{n+1},\label{405b}
\end{eqnarray}
where $\beta$ and $\varsigma$ are parameters of the methods that will be fixed later.
Approximating 
${\mathbf{G}}_1({\varphi}_h)$ and ${\mathbf{G}}_2({\psi}_h)$
by 
$\widetilde{\mathbf{G}}_1({\varphi}_h^n)=\sum\limits_{j=0}^N \delta t g_1^j\mathbf{D}^2_0\varphi^{n-j}$
and
$\widetilde{\mathbf{G}}_2({\psi}_h^n)=\sum\limits_{j=0}^N \delta t g_2^j\mathbf{D}^+\mathbf{D}^-\psi^{n-j}$ 
(for $N$ large enough), and
replacing \eqref{402}-\eqref{405b} in the equation of motion \eqref{401}, we obtain
\begin{eqnarray}
\left(\mathbf{M}+\varsigma\delta t\,\mathbf{C}+\beta\delta t^2
(\mathbf{K}+\mathbf{G}_0)\right)
\left[\begin{array}{c}
\ddot{\varphi}_h^{n+1}\\ \ddot{\psi}_h^{n+1} \\ \ddot{v}_h^{n+1}
\end{array}\right]
=
-\mathbf{C}\left( \left[\begin{array}{c}
\dot{\varphi}_h^{n}\\ \dot{\psi}_h^{n} \\ \dot{v}_h^{n}
\end{array}\right] + (1-\varsigma)\delta t\,\left[\begin{array}{c}
\ddot{\varphi}_h^{n}\\ \ddot{\psi}_h^{n} \\ \ddot{v}_h^{n}
\end{array}\right]\right)\nonumber
\\
- (\mathbf{K}+\mathbf{G}_0) \left(
\left[\begin{array}{c}
{\varphi}_h^{n}\\ {\psi}_h^{n} \\ {v}_h^{n}
\end{array}\right]
+\delta t \left[\begin{array}{c}
\dot{\varphi}_h^{n}\\ \dot{\psi}_h^{n} \\ \dot{v}_h^{n}
\end{array}\right] +
\left(\frac{1}{2}-\beta\right)\delta t^2 \left[\begin{array}{c}
\ddot{\varphi}_h^{n}\\ \ddot{\psi}_h^{n} \\ \ddot{v}_h^{n}
\end{array}\right]
\right)\label{406}
\\
-\sum_{j=1}^{N}\mathbf{G}_j
\left[\begin{array}{c}
{\varphi}_h^{n+1-j}\\ {\psi}_h^{n+1-j} \\ {v}_h^{n+1-j}
\end{array}\right],\nonumber
\end{eqnarray}
where
$\mathbf{G}_j=\begin{bmatrix}
\delta t g_1^j\mathbf{D}^2_0 & 0 & 0 \\ 0 & \delta t g_2^j\mathbf{D}^+\mathbf{D}^- & 0 \\ 0 & 0 & 0 
\end{bmatrix} $, $j=0,\ldots,N$.
We introduce the variables 
\begin{equation}
\begin{split}
\eta^{n,j}:=\varphi^n-\varphi^{n-j}\\
z^{n,j}:=\psi^n-\psi^{n-j}
\end{split}
\label{eq410}
\end{equation}
which verify
\begin{equation}
\begin{split}
\eta^{n,j}-\eta^{n,j-1}=\varphi^{n+1-j}-\varphi^{n-j},\\
z^{n,j}-z^{n,j-1}=\psi^{n+1-j}-\psi^{n-j}.
\end{split}
\label{eq411}
\end{equation}
Thus, from \eqref{402} and \eqref{405}, we obtain
\begin{equation}
\begin{split}
	\delta t\dot{\varphi}^{n+\frac{1}{2}}=\eta^{n+1,j}-\eta^{n,j-1}+
	\delta t^2\left(\beta-\frac{1}{2}\varsigma\right)\left(\ddot{\varphi}^{n+1}-\ddot{\varphi}^{n}\right)\\
		\delta t\dot{\psi}^{n+\frac{1}{2}}=z^{n+1,j}-z^{n,j-1}+
		\delta t^2\left(\beta-\frac{1}{2}\varsigma\right)\left(\ddot{\psi}^{n+1}-\ddot{\psi}^{n}\right)\\
\end{split}
\label{eq412}
\end{equation}
where $\vartheta^{n+\frac{1}{2}}:=\frac{\vartheta^n+\vartheta^{n+1}}{2}$,
for all $\vartheta^n$, with $n\in\mathbb{Z}$.
We will need the following Lemma:
\vskip0,1truecm
\begin{lemma}
	Let $\beta=\frac{1}{2}\varsigma$. Then
	\begin{eqnarray}
&&		\begin{split}
\sum_{j=1}^N\chi^j\mathbf{D}^2_0\varphi^{n+1-j}\cdot\dot{\varphi}^{n+\frac{1}{2}} =&
\frac{1}{2\delta t}	\Bigg[-\left(\sum_{j=1}^{N}\chi^j\right)\left(\|\mathbf{D}^+\varphi\|^2  + \left(\frac{\varphi_1}{h^2}\right)^2\right)
	\\ &
	+\sum_{j=1}^N\chi^j\left(\|\mathbf{D}^+\eta^{\cdot,j}\|^2+\left(\frac{\eta^{\cdot,j}_1}{h^2}\right)^2\right)\Bigg]_{n}^{n+1}
	\\
& - \frac{1}{2\delta t}\sum_{j=1}^{N}\left(\chi^{j+1}-\chi^j\right)\left(\| \mathbf{D}^+\eta^{n,j}\|^2 + \left(\frac{\eta^{n,j}_1}{h^2}\right)^2\right)\\
	&
	+\frac{1}{2\delta t}\chi^{N+1}\left(\| \mathbf{D}^+\eta^{n,N}\|^2 + \left(\frac{\eta^{n,N}_1}{h^2}\right)^2\right),
	\end{split}
	\label{eq413}
	\\
&&	\begin{split}
\sum_{j=1}^N\chi^j\mathbf{D}^+\mathbf{D}^-\psi^{n+1-j}\cdot\dot{\psi}^{n+\frac{1}{2}} =&
\frac{1}{2\delta t}\left[-\left(\sum_{j=1}^{N}\chi^j\right)\|\mathbf{D}^-\psi\|^2
+\sum_{j=1}^N\chi^j\|\mathbf{D}^-z^{\cdot,j}\|^2\right]_{n}^{n+1}
\\
& - \frac{1}{2\delta t}\sum_{j=1}^{N}\left(\chi^{j+1}-\chi^j\right)\| \mathbf{D}^-z^{n,j}\|^2
+ \frac{1}{2\delta t}\chi^{N+1}\| \mathbf{D}^-z^{n,N}\|^2,
	\end{split}
	\label{eq414}
	\end{eqnarray}
	for all $\left(\chi_j\right)_{j\in \mathbb{N}}$
	with $\chi_j\in\mathbb{R}$.
\end{lemma}
\vskip0,1truecm
\begin{proof}
Let $I(\psi^{n}):=\sum\limits_{j=1}^N\chi^j\mathbf{D}^+\mathbf{D}^-\psi^{n-j}\cdot\dot{\psi}^{n+\frac{1}{2}}$,
	and
$I(\psi^{n+\frac{1}{2}}):=\sum\limits_{j=1}^N\chi^j\mathbf{D}^+\mathbf{D}^-\psi^{n+\frac{1}{2}-j}\cdot\dot{\psi}^{n+\frac{1}{2}}$.
Then, using \eqref{eq410}, \eqref{eq411} and \eqref{eq412}, with $\beta=\frac{1}{2}\varsigma$,
we obtain
\begin{eqnarray}
\lefteqn{
I(\psi^{n+\frac{1}{2}})=
-\sum\limits_{j=1}^N\chi^j \mathbf{D}^-\left( \psi^{n+\frac{1}{2}}-
z^{j,n+\frac{1}{2}}\right)\cdot
 \mathbf{D}^- \frac{\psi^{n+1}-\psi^n}{\delta t}
}\nonumber\\
&=&
-
\overbrace{
\frac{1}{2\delta t}\left(\sum_{j=1}^{N}\chi^j\right)\left(\|\mathbf{D}^-\psi^{n+1}\|^2
- \|\mathbf{D}^-\psi^{n}\|^2\right)}^{A}
+ 
\overbrace{
\frac{1}{\delta t} \sum_{j=1}^{N}\chi^j\mathbf{D}^-z^{n+\frac{1}{2}}\cdot
\mathbf{D}^-\left(z^{n+1,j}-z^{n,j}\right)
}^{B} \nonumber \\
&&
+ 
\underbrace{
 \frac{1}{\delta t} \sum_{j=1}^{N}\chi^j\mathbf{D}^-z^{n+\frac{1}{2}}\cdot
\mathbf{D}^-\left(z^{n,j}-z^{n,j-1}\right). 
}_{C}
\label{eq415}
\end{eqnarray}
The term $-A+B$ corresponds to the right hand side of the term in bracket $\Big[\cdot\Big]^{n+1}_n$
in \eqref{eq414}. On the other hand, the term $C$ can be written as
\begin{eqnarray*}
C&=&
 \frac{1}{2\delta t} \sum_{j=1}^{N}\chi^j
 \mathbf{D}^-\left(z^{n,j}+z^{n,j-1}\right)\cdot
 \mathbf{D}^-\left(z^{n,j}-z^{n,j-1}\right)\\ &&
 +
  \frac{1}{2\delta t} \sum_{j=1}^{N}\chi^j
  \mathbf{D}^-\left(z^{n+1,j}-z^{n,j-1}\right)\cdot
  \mathbf{D}^-\left(z^{n,j}-z^{n,j-1}\right)\\
  &=&
  - \frac{1}{2\delta t}\sum_{j=1}^{N}\left(\chi^{j+1}-\chi^j\right)\| \mathbf{D}^-z^{n,j}\|^2
  + \frac{1}{2\delta t}\chi^{N+1}\| \mathbf{D}^-z^{n,N}\|^2
  + I(\psi^n) - I(\psi^{n+\frac{1}{2}})
\end{eqnarray*}
Thus, observing that
the left hand side of \eqref{eq414} is given by
$I(\psi^{n+1})=2I(\psi^{n+\frac{1}{2}})- I(\psi^n)$, and
replacing this last expression in \eqref{eq415}, we obtain \eqref{eq414}. 
Repeating the same calculations for the couple$(\varphi,\eta)$, and
	considering \eqref{eq43b}, it follows  \eqref{eq415}.
\end{proof}
\vskip0,1truecm
The acceleration $[\ddot{\varphi}_{h}^{n+1}, \ddot{\psi}_{h}^{n+1}, \ddot{v}_{h}^{n+1} ]^\top$
is computed from \eqref{406}, and the velocities
$[\dot{\varphi}_{h}^{n+1}, \dot{\psi}_{h}^{n+1}, \dot{v}_{h}^{n+1} ]^\top$
are obtained from \eqref{402} and \eqref{404}, respectively.
Finally, the displacement
$[{\varphi}_{h}^{n+1}, {\psi}_{h}^{n+1}, {v}_{h}^{n+1} ]^\top$
follows from \eqref{403} and \eqref{405},
by simple matrix operations.
Thus, the fully discrete energy of the system \eqref{402}-\eqref{406} is given by
\begin{eqnarray}
\label{Eh}
\mathcal{E}_h^{n} &:=&
\frac{1}{2} [\dot{\varphi}_{h}^{n}, \dot{\psi}_{h}^{n}, \dot{v}_{h}^{n} ]^\top
\mathbf{M}
\left[\begin{array}{c}
\dot{\varphi}_h^{n}\\ \dot{\psi}_h^{n} \\ \dot{v}_h^{n}
\end{array}\right]
+
\frac{1}{2} [{\varphi}_{h}^{n}, {\psi}_{h}^{n}, {v}_{h}^{n} ]^\top
\mathbf{K}
\left[\begin{array}{c}
{\varphi}_h^{n}\\ {\psi}_h^{n} \\ {v}_h^{n}
\end{array}\right]
+\mathcal{G}\left(\left[\begin{array}{c}
{\varphi}_h^{n}\\ {\psi}_h^{n} \\ {v}_h^{n}
\end{array}\right],\left[\begin{array}{c}
{\varphi}_h^{n}\\ {\psi}_h^{n} \\ {v}_h^{n}
\end{array}\right] \right)
\\
\nonumber
&=&
\frac{3}{2}\rho_1 \|\dot{\varphi}_h^{n}\|^2 +
\frac{3}{2}\rho_2 \|\dot{\psi}_h^{n}\|^2 +
\frac{3}{2} b  \| \mathbf{D}^-{\varphi}_h^{n}\|^2 +
\frac{1}{2}\rho_2 \|\dot{v}_h^{n}\|^2 +
\frac{3}{2} b  \| \mathbf{D}^-v_h^{n}\|^2 \\
&& +
\frac{3}{2} k  \|\mathbf{P}^- v_h^{n} - \mathbf{P}^- \psi_h^{n} - \mathbf{D}^+\varphi_h^{n} \|^2 
+  \frac{3\,k}{2\,h^2} \varphi_1^2
\nonumber
\\ && +
\frac{1}{2}	\Bigg[\left(\frac{1}{2}g^0_1-\sum_{j=1}^{N}g_1^{j-\frac{1}{2}}\right)\left(\|\mathbf{D}^+\varphi^n_h\|^2  + \left(\frac{\varphi^n_{h,1}}{h^2}\right)^2\right)
	\nonumber	\\ &&
	+\sum_{j=1}^N g_1^{j-\frac{1}{2}}\left(\|\mathbf{D}^+\eta^{n,j}_h\|^2+\left(\frac{\eta^{n,j}_{h,1}}{h^2}\right)^2\right)\Bigg]
	\nonumber	\\ && +
	\frac{1}{2}\left[\left(\frac{1}{2}g^0_2-\sum_{j=1}^{N}g_2^{j-\frac{1}{2}}\right)\|\mathbf{D}^-\psi^n_h\|^2
	+\sum_{j=1}^Ng_2^{j-\frac{1}{2}}\|\mathbf{D}^-z^{n,j}_h\|^2\right]
	\nonumber
\end{eqnarray}
where $\mathcal{G}(\cdot,\cdot)$ is the bilinear form derived from the memory term and described in the following line, and $g_{1,2}^{j-{\frac{1}{2}}}=\frac{1}{2}(g_{1,2}^{j-1}+g_{1,2}^{j})$.
This is an approximation of energy for the continuous case.
The increment of this energy can be expressed in terms of mean values and increments of the displacement and velocity. 
Then, we choose ${\varsigma=\frac{1}{2}}$ and ${\beta=\frac{\varsigma}{2}}$, reducing
the above expression to
\begin{eqnarray*}
\lefteqn{	\mathcal{E}_\delta^{n+1}-\mathcal{E}_\delta^n =
	-4\,\delta t \,\gamma\,\|\dot{v}^{n+\frac{1}{2}}_h\|^2}	\\
	&&
\frac{1}{4}\sum_{j=1}^{N}\left(g_1^{j+1}-g_1^{j-1}\right)\left(\| \mathbf{D}^+\eta^{n,j}\|^2 + \left(\frac{\eta^{n,j}_1}{h^2}\right)^2\right)\\
	&&
	-\frac{1}{4}g_1^{0}\left(\| \mathbf{D}^+\eta^{n+1,1}\|^2 + \left(\frac{\eta^{n+1,1}_1}{h^2}\right)^2\right) 
	-\frac{1}{2}g_1^{N+\frac{1}{2}}\left(\| \mathbf{D}^+\eta^{n,N}\|^2 + \left(\frac{\eta^{n,N}_1}{h^2}\right)^2\right) 
	\\ &&
+ \frac{1}{4}\sum_{j=1}^{N}\left(g_2^{j+1}-g_2^{j-1}\right)\| \mathbf{D}^-z^{n,j}\|^2
- \frac{1}{4}g_2^{0}\| \mathbf{D}^-z^{n+1,1}\|^2
- \frac{1}{2}g_2^{N+\frac{1}{2}}\| \mathbf{D}^-z^{n,N}\|^2
	\ {\leqslant}\ {0}.
	\label{407}
\end{eqnarray*}
\vskip0,1truecm
With this, the fully discrete energy obtained by the $\beta-$Newmark method is decreasing and we
expect that its asymptotic behavior be a reflection of the continuous case (see \cite{Krenk}
and also \cite{arsv1,arsv2}).

\begin{figure}[h]
	\begin{tabular}{ll}\hspace{-0.5cm}
		\includegraphics[scale=0.55]{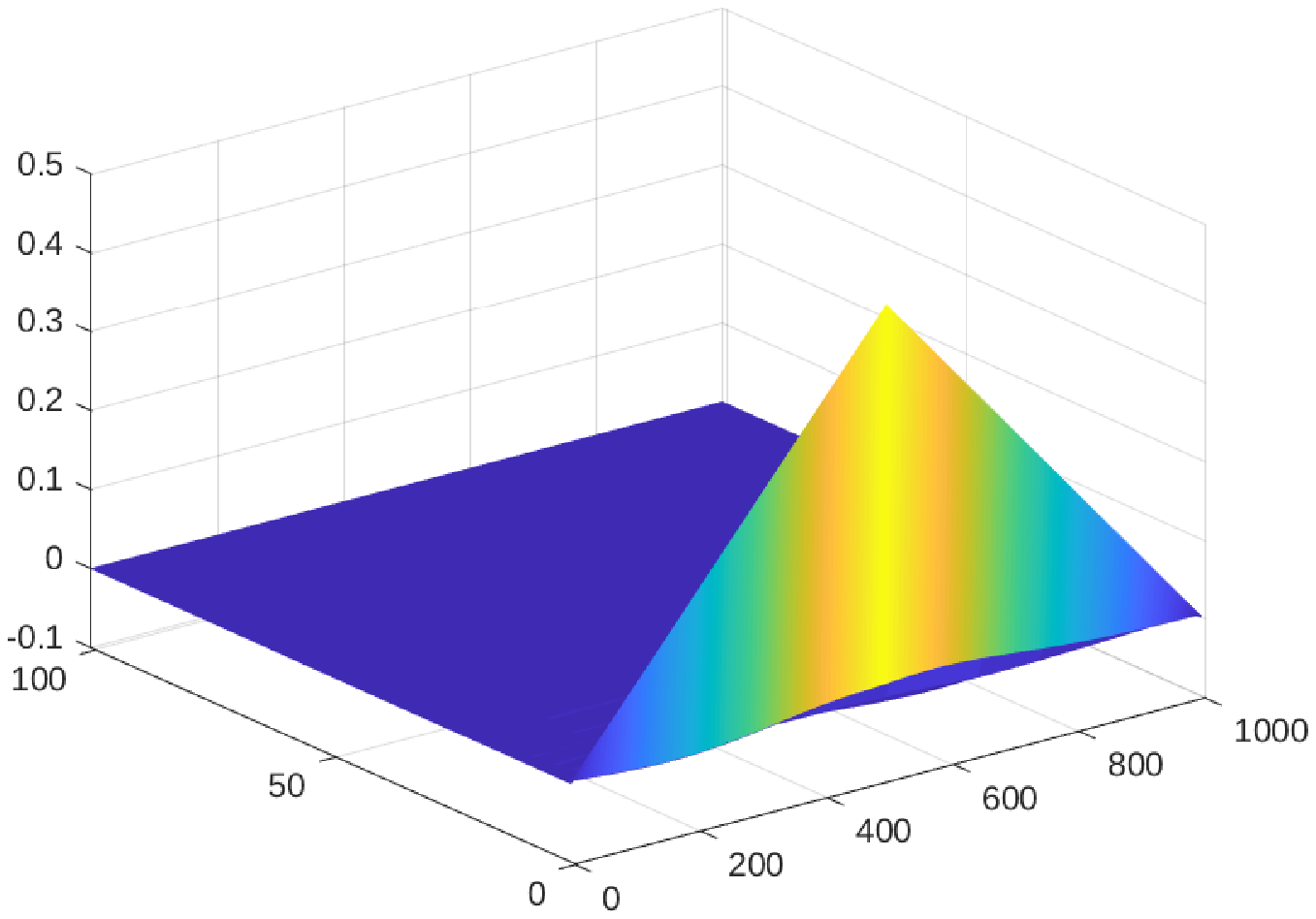} & \hspace{-0.5cm}
		\includegraphics[scale=0.55]{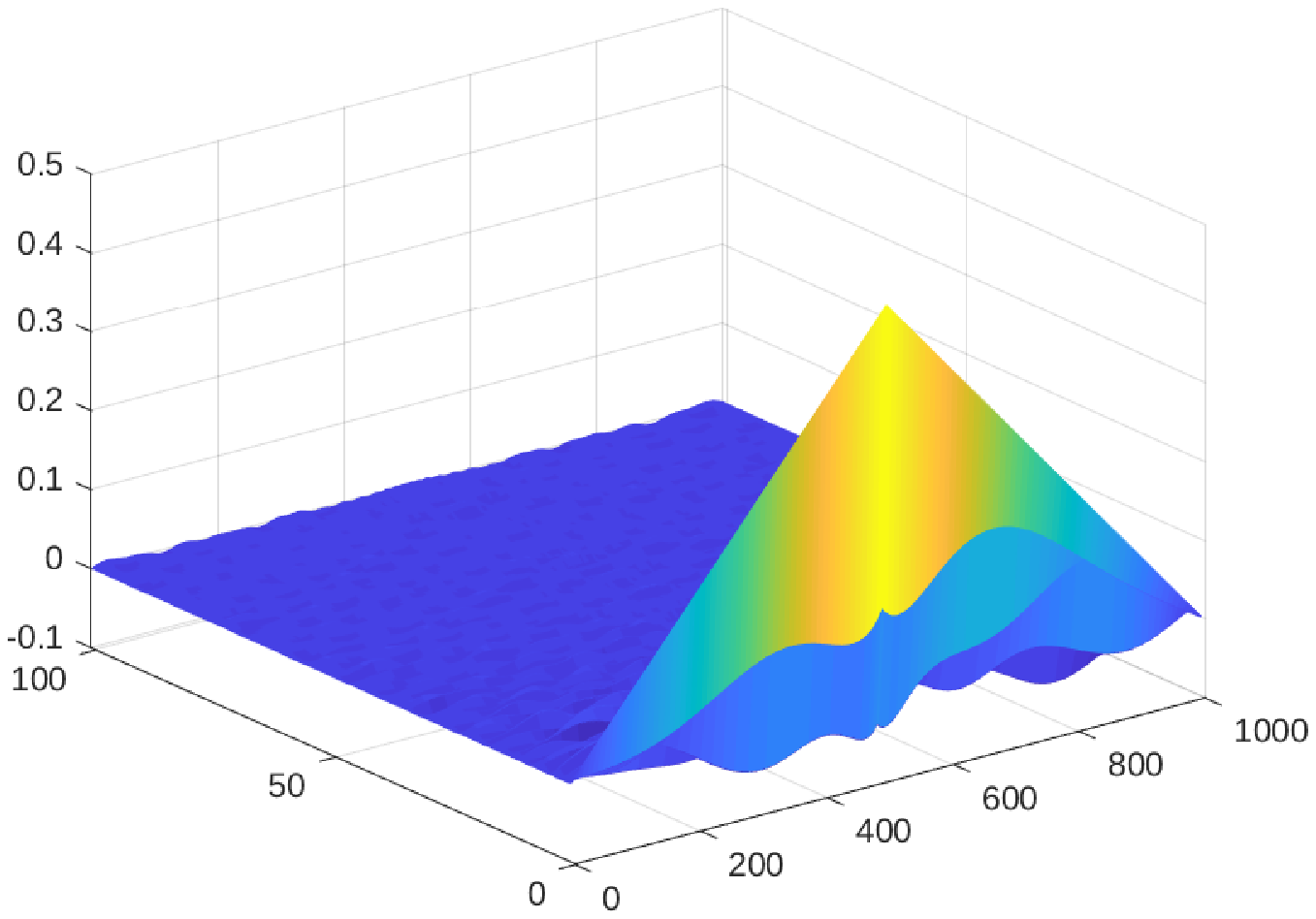} \\ \hspace{-0.5cm}
	\end{tabular}
	\caption{ $\varphi(x,t)$ and $\varphi_t(x,t)$ \label{Fig1}}
\end{figure}

\begin{figure}[h]
	\begin{tabular}{ll}\hspace{-0.5cm}
		\includegraphics[scale=0.55]{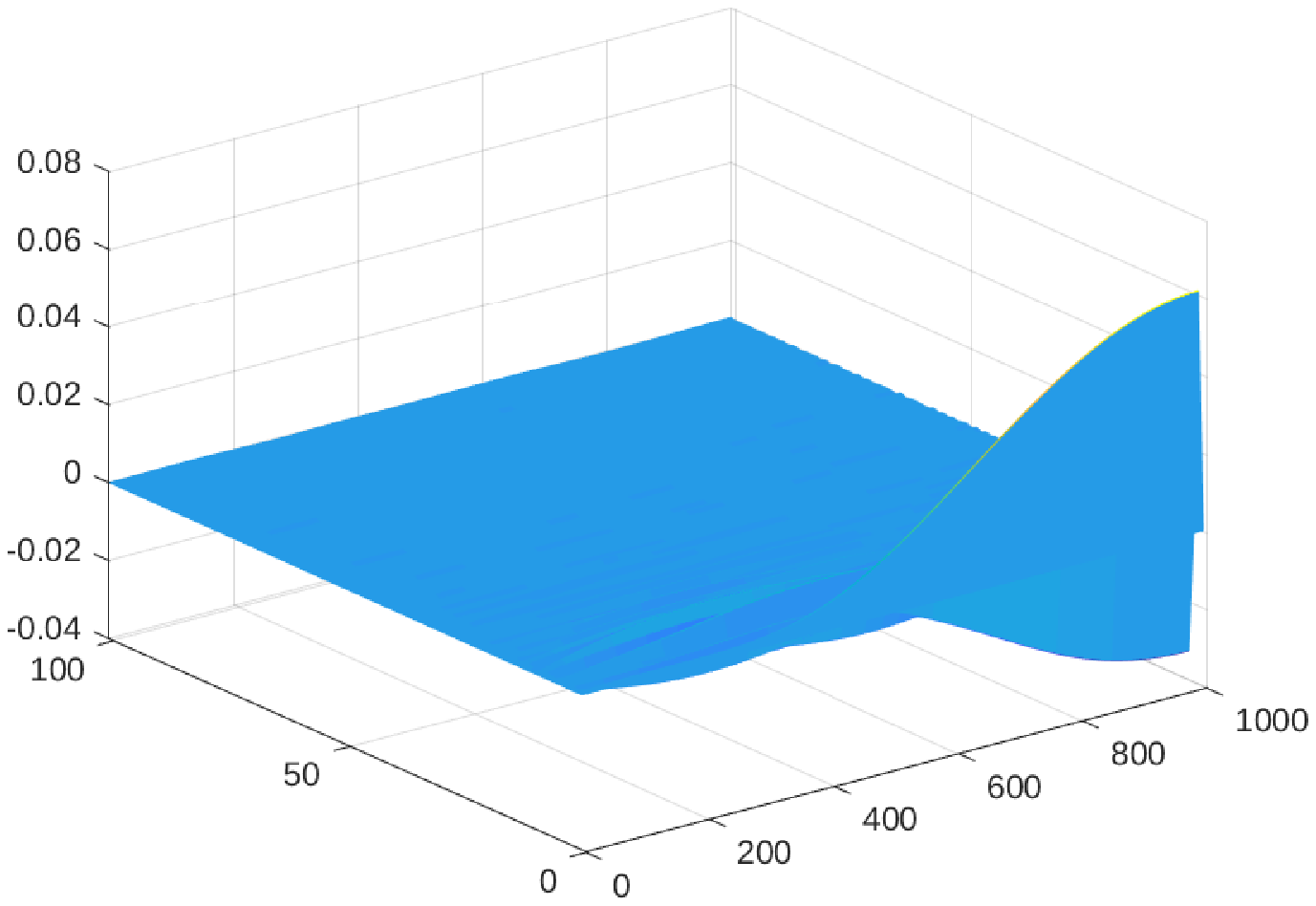} & \hspace{-0.5cm}
		\includegraphics[scale=0.55]{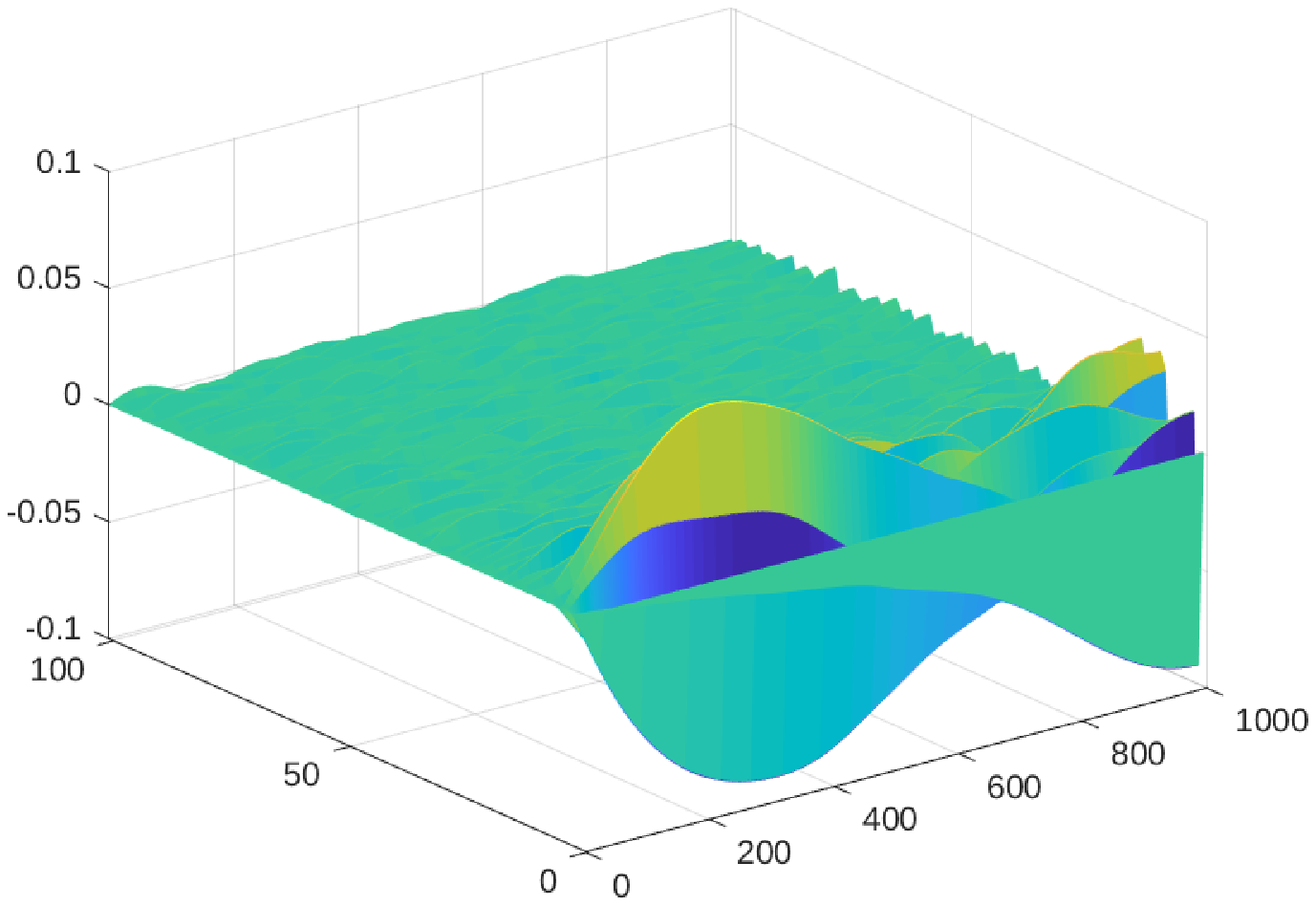} \\ \hspace{-0.5cm}
	\end{tabular}
	\caption{ $\psi(x,t)$ and $\psi_t(x,t)$ \label{Fig2}}
\end{figure}

\begin{figure}[h]
	\begin{tabular}{ll}\hspace{-0.5cm}
		\includegraphics[scale=0.55]{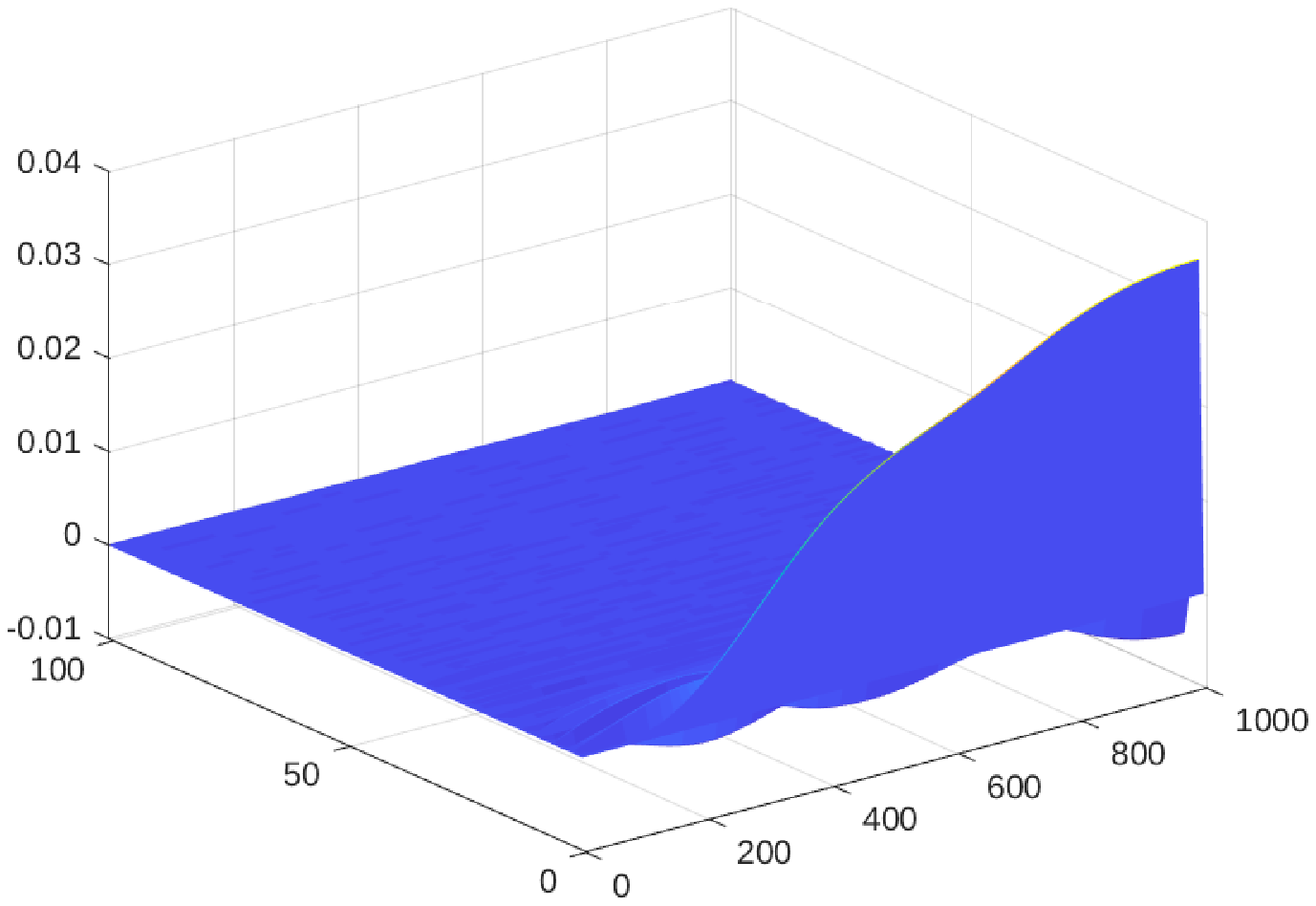} & \hspace{-0.5cm}
		\includegraphics[scale=0.55]{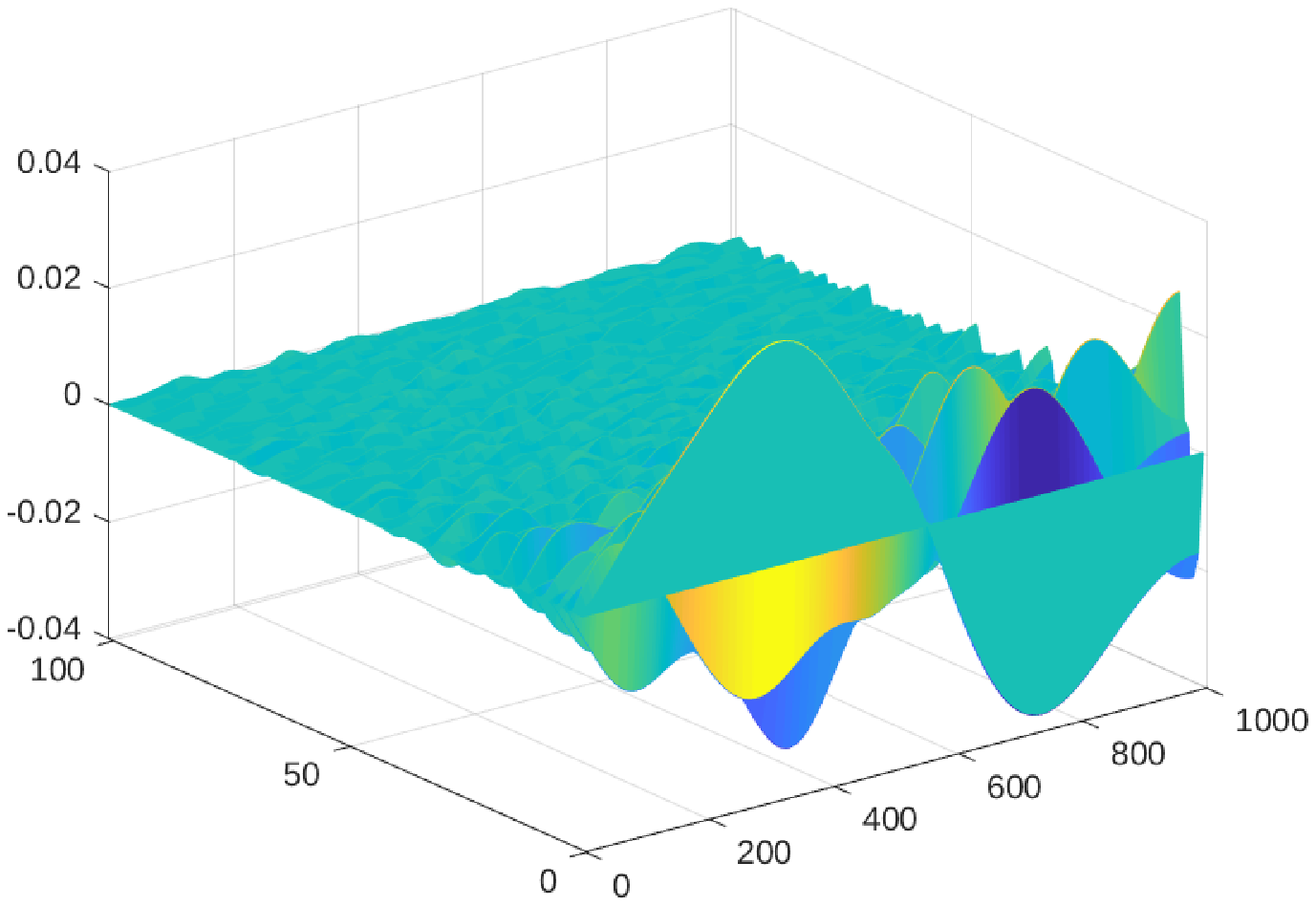} \\ \hspace{-0.5cm}
	\end{tabular}
	\caption{ $v(x,t)$ and $v_t(x,t)$ \label{Fig3}}
\end{figure}

\begin{figure}[h]
		\includegraphics[scale=.5]{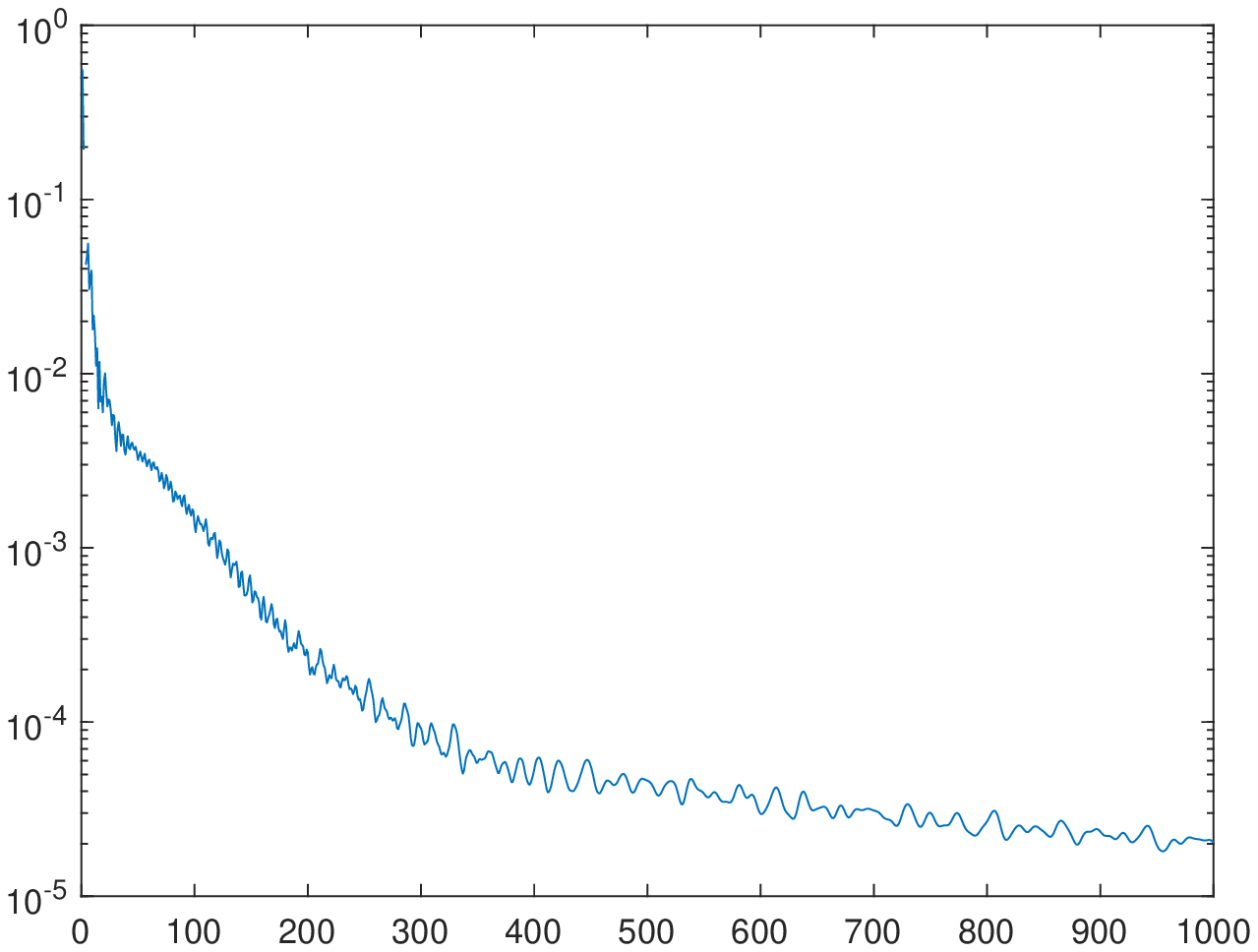}
	\caption{
		Decay Energy
		\label{Fig4}}
\end{figure}

\section{Comments and issues}

{\bf Comment 1.} The speeds of wave propagations of the both last two equations in \eqref{109} are equal to $\sqrt{\frac{b}{\rho_2}}$.
Our results hold true when the last two equations in \eqref{109} have different speeds of wave propagations (that is when $b$ in the last equation in \eqref{109} is replaced by ${\tilde b}>0$).
\vskip0,1truecm
{\bf Comment 2.} Our results hold true when $\delta =0$.
\vskip0,1truecm 
{\bf Comment 3.} The last equation in \eqref{109} can be controlled via an infinite memory 
\begin{equation*}
\int_{0}^{\infty} g_3 (s)\,v_{xx} (t - s)\ ds
\end{equation*}
instead of the linear frictional damping $4\gamma v_t$, where $g_3\,: \mathbb{R}_{+}\to \mathbb{R}_{+}$ is a given relaxation function satisfying the same hypotheses as $g_1$ and $g_2$. To prove the well-posedness results, we introduce a third variable $w$ similar to $\eta$ and $z$ given by 
\begin{equation*}
w(x,\,t,\,s) = v(x,\,t) - v(x,\,t - s),
\end{equation*}
we define its space $L_{g_3}$ as $L_{g_2}$ and we do some logical modifications. For the stability result, we add to $E$ in \eqref{214} 
the term
\begin{equation*}
\frac{1}{2} \int_{0}^{\infty}g_3 (s)\Vert w_x\Vert^2\, ds
\end{equation*} 
and we add to the definition of $F$ in \eqref{402} the integral 
\begin{equation*}
-\rho_2\int_{0}^1 v_{t}\int_{0}^{\infty}g_3 (s)\,w\, ds\,dx.
\end{equation*}  
\vskip0,1truecm
{\bf Acknowledgment.} This work was initiated during the visit in July-August 2017 of the first author to LNCC and RJ university, Brazil, and B\'{\i}o-B\'{\i}o and Concepci\'{o}n universities, Chile, and finished during the visit in June 2018 of the third author to Lorraine-Metz university, France, and the visit in August 2018 of the first author to B\'{\i}o-B\'{\i}o and Concepci\'{o}n universities, Chile. The first and third authors thank LNCC, Lorraine-Metz, RJ, B\'{\i}o-B\'{\i}o and Concepci\'{o}n universities for their kind support and hospitality. This work was supported FONDECYT grant no. 1180868, and by ANID-Chile through the project {\sc Centro de Modelamiento Matem\'atico} (AFB170001)  of the PIA Program: Concurso Apoyo a Centros Científicos y Tecnológicos de Excelencia con Financiamiento Basal.

\end{document}